\documentclass{iopart}

\usepackage[latin1]{inputenc}
\usepackage{iopams}
\usepackage{harvard}
\usepackage{graphicx}
\usepackage{psfrag}

\newtheorem{theorem}{Theorem}[section]
\newtheorem{mainthm}{Main Theorem}
\newtheorem{lemma}[theorem]{Lemma}

\newtheorem{corollary}[theorem]{Corollary}
\newtheorem{remark}[theorem]{Remark}

\newtheorem{definition}[theorem]{Definition}

\newcommand{\field}[1]{\mathbb{#1}}

\def\P{\tilde{P}}
\def\C{\mathbb{C}}

\def\eps{\epsilon}

\def\fC{\field{C}}

\def\fN{\field{N}}

\def\fR{\field{R}}

\def\fZ{\field{Z}}

\def\cA{\mathcal{A}}
\def\cB{\mathcal{B}}
\def\cC{\mathcal{C}}
\def\cD{\mathcal{D}}

\def\cF{\mathcal{F}}

\def\cK{\mathcal{K}}
\def\cI{\mathcal{I}}

\def\cM{\mathcal{M}}
\def\cT{\mathcal{T}}
\def\cR{\mathcal{R}}

\def\cP{\mathcal{P}}
\def\cS{\mathcal{S}}
\def\cU{\mathcal{U}}
\def\cW{\mathcal{W}}
\def\cV{\mathcal{V}}
\def\cX{\mathcal{X}}
\def\cZ{\mathcal{Z}}
\def\cO{\mathcal{O}}

\def\bLambda{{\bf \Lambda_*}}
\def\blambda{{\boldsymbol \lambda_*}}
\def\bmu{{\boldsymbol \mu_*}}
\def\bF{{\bf F_*}}
\def\bG{{\bf G_*}}

\def\tbG{{\bf \tilde{G}_*}}
\def\tG{{\tilde{G}}}
\def\bs{{\bf s^*}}
\def\bA{{\bf A}}

\def\bW{{\bf W}}
\def\bX{{\bf X}}
\def\bY{{\bf Y}}

\def\bk{{\bf k}}
\def\bp{{\bf p}}
\def\bq{{\bf q}}
\def\tp{{\tilde p}}
\def\tD{{\tilde \Delta}}
\def\tF{{\tilde F}}
\def\loc{{\rm loc}}
\def\hF{{\hat F}}
\newcommand{\cover}[1]{\stackrel{#1}{\Longrightarrow}} 
\newcommand{\bcover}[1]{\stackrel{#1}{\Longleftarrow}} 
\newcommand{\gcover}[1]{\stackrel{#1}{\Longleftrightarrow}}
\newcommand{\converge}[1]{\stackrel{#1}{\longrightarrow}}

\newcommand{\Bm}[2]{\!\!\!\phantom{a}_{#1}B_1^{#2}}
\newcommand{\BmZ}[2]{\!\!\!\phantom{a}_{#1}B_0^{#2}}

\date{2009-03-31}

\begin{document}
\title[Dynamics of the Universal Area-Preserving Map: Hyperbolic Sets]{Dynamics of the Universal Area-Preserving Map Associated with Period Doubling: Hyperbolic Sets}
\author{Denis Gaidashev, Tomas Johnson}
\address{Department of Mathematics, Uppsala University, Box 480, 751 06 Uppsala, Sweden}
\eads{\mailto{gaidash@math.uu.se}, \mailto{tomas.johnson@math.uu.se}}
\begin{abstract}
It is known that the famous Feigenbaum-Coullet-Tresser period doubling universality has a counterpart for area-preserving maps of ${\fR}^2$. A renormalization approach has been used in \cite{EKW1} and \cite{EKW2} in a computer-assisted proof of existence of a ``universal'' area-preserving map $F_*$ ---  a map with orbits of all binary periods $2^k, k \in \fN$.  In this paper, we consider maps in some neighbourhood of  $F_*$ and study their dynamics. 

We first demonstrate that the map $F_*$ admits a ``bi-infinite heteroclinic tangle'': a sequence of periodic points $\{z_k\}$, $k \in \fZ$, 
\begin{equation} 
|z_k|   \converge{{k \rightarrow \infty}}  0, \quad  |z_k|    \converge{{ k \rightarrow -\infty}} \infty,
\end{equation}
 whose stable and unstable manifolds intersect transversally; and, for any $N \in \fN$, a compact invariant set on which $F_*$ is homeomorphic to a topological Markov chain on  the space of all two-sided sequences composed of $N$ symbols.
A corollary of these results is the existence of {\it unbounded} and {\it oscillating} orbits.

We also show that the third iterate for all maps close to $F_*$ admits a horseshoe. We use distortion tools to provide rigorous bounds on the Hausdorff dimension of the associated locally maximal invariant hyperbolic set:
$$ 0.7673 \ge  {\rm dim}_H(\cC_F) \ge \varepsilon \approx 0.00013 \, e^{-7499}.$$ 

\end{abstract}
\submitto{Nonlinearity}

\maketitle

\setcounter{page}{1}


\section{Introduction}

Following the pioneering discovery of the Feigenbaum-Coullet-Tresser period doubling universality in unimodal maps \cite{Fei1}, \cite{Fei2}, \cite{TC} universality has been demonstrated to be a rather generic phenomenon in dynamics.

To prove universality one usually introduces a {\it renormalization} operator on a functional space, and demonstrates that this operator has a hyperbolic fixed point.

Such renormalization approach to universality has been very successful in one-dimensional dynamics, and has led to explanation of universality in unimodal maps \cite{Eps1}, \cite{Eps2},\cite{Lyu}, critical circle maps \cite{dF1,dF2}, \cite{Ya1}, \cite{Ya2} and holomorphic maps with a  Siegel disk \cite{McM}, \cite{Ya3}, \cite{GaiYa}.

Universality has been abundantly observed in higher dimensions, in particular, in two and more dimensional dissipative systems (cf.  \cite{CEK1}, \cite{Spa}), in area-preserving maps, both as the period-doubling universality \cite{DP}, \cite{Hel}, \cite{BCGG}, \cite{CEK2}, \cite{EKW1}, \cite{EKW2}, \cite{GK1} and as the universality associated with the break-up of invariant surfaces  \cite{Shen}, \cite{McK1}, \cite{McK2}, \cite{ME},  and in Hamiltonian flows \cite{ED},\cite{AK}, \cite{AKW}, \cite{Koch1}, \cite{Koch2}, \cite{Koch3}, \cite{GK}, \cite{Gai1}, \cite{Kocic} .

It has been established that the universal behaviour in dissipative and conservative higher dimensional systems is fundamentally different. The case of the dissipative systems is often reducible to the one-dimensional Feigenbaum-Coullet-Tresser universality (\cite{CEK1}, \cite{dCLM}, \cite{LM}). Universality for highly dissipative H\'enon like maps
$$F(x,u)=(f(x)-\epsilon(x,u),x),$$
where $f$ is a unimodal map, $\epsilon$ - sufficiently small together with its derivatives, has been demonstrated in \cite{dCLM}. Specifically, it has been shown that these maps are in the Feigenbaum-Coullet-Tresser universality class, and that the degenerate map
$$F^*(x,u)=(f^*(x),x),$$ 
where $f^*$ is the  Feigenbaum-Coullet-Tresser fixed point, is a renormalization fixed point in that class. The authors of \cite{dCLM} have also constructed a Cantor set (an attractor) for infinitely renormalizable maps, and  demonstrated that for such H\'enon-like maps universality coexists with non-rigidity: the H\"older exponent of the conjugacy between the actions of any two infinitely renormalizable maps with nonequal average Jacobians on their Cantor sets has an upper bound less than $1$.

The case of area-preserving maps seems to be very different, and at present there is no deep understanding of universality in conservative systems, other than in the ``trivial'' case of the universality for systems ``near integrability'' \cite{Koch1}, \cite{Koch2}, \cite{Gai1}, \cite{Kocic}, \cite{KLDM}.

An infinite period-doubling cascade in families of area-preserving maps was observed by several authors in  early 80's \cite{DP}, \cite{Hel}, \cite{BCGG}, \cite{Bou}, \cite{CEK2}. The period-doubling phenomenon can be illustrated with the area-preserving H\' enon family (cf. \cite{Bou}) :
$$ H_a(x,u)=(-u +1 - a x^2, x).$$

Maps $H_a$ have a fixed point $((-1+\sqrt{1+a})/a,(-1+\sqrt{1+a})/a) $ which is stable for $-1 < a < 3$. When $a_1=3$ this fixed point becomes unstable, at the same time an orbit of period two is born with $H_a(x_\pm,x_\mp)=(x_\mp,x_\pm)$, $x_\pm= (1\pm \sqrt{a-3})/a$. This orbit, in turn, becomes unstable at $a_2=4$, giving birth to a period $4$ stable orbit. Generally, there  exists a sequence of parameter values $a_k$, at which the orbit of period $2^{k-1}$ turns unstable, while at the same time a stable orbit of period $2^k$ is born. The parameter values $a_k$ accumulate on some $a_\infty$. The crucial observation is that the accumulation rate
\begin{equation}
\lim_{k \rightarrow \infty}{a_k-a_{k-1} \over  a_{k+1}-a_k } = 8.721...
\end{equation} 
is universal for a large class of families, not necessarily H\'enon.

Furthermore, the $2^k$ periodic orbits scale asymptotically with two scaling parameters
\begin{equation}
\lambda=-0.249 \ldots,\quad \mu=0.061 \ldots
\end{equation}

\begin{figure}[ht]
\begin{center}
\includegraphics[angle=-90,width=0.7 \textwidth]{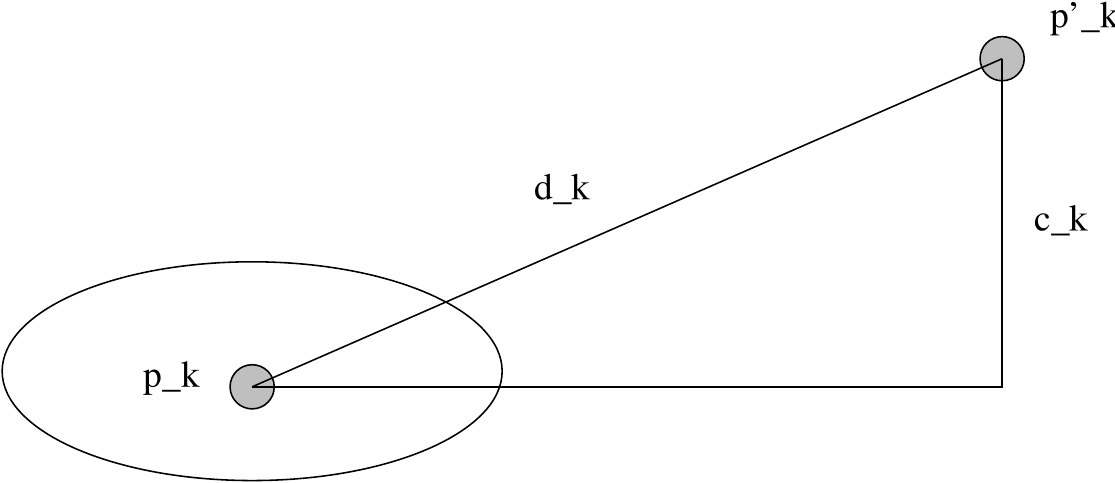}
\caption{The geometry of the period doubling. $p_k$ is the further elliptic point that has bifurcated from the hyperbolic point $p'_k$.}\label{bifGeom}
\end{center}
\end{figure}

To explain how orbits scale with $\lambda$ and $\mu$ we will follow \cite{Bou}. Consider an interval $(a_k,a_{k+1})$ of parameter values in a ``typical'' family $F_a$. For any value $\alpha \in (a_k,a_{k+1})$ the map $F_\alpha$ possesses a stable periodic orbit of period $2^{k}$. We fix some $\alpha_k$ within the interval $(a_k,a_{k+1})$ in some consistent way; for instance, by requiring that the restriction of $F^{2^{k}}_{\alpha_k}$ to a neighbourhood of a stable periodic point in the $2^{k}$-periodic orbit is conjugate, via a diffeomorphism $H_k$, to a rotation with some fixed rotation number $r$.  Let $p'_k$ be some unstable periodic point in the $2^{k-1}$-periodic orbit, and let $p_k$ be the further of the two stable $2^{k}$-periodic points that bifurcated from $p'_k$.  Denote with $d_k=|p'_k-p_k|$, the distance between $p_k$ and $p'_k$. The new elliptic point $p_k$ is surrounded by invariant ellipses; let $c_k$ be the distance between $p_k$ and $p'_k$ in the direction of
  the minor semi-axis of an invariant ellipse surrounding $p_k$, see Figure \ref{bifGeom}. Then,
$$
{1 \over \lambda}=-\lim_{k \rightarrow \infty}{ d_k \over d_{k+1}},\quad 
{\lambda \over \mu}=-\lim_{k \rightarrow \infty}{ \rho_k \over \rho_{k+1}}, \quad 
{1 \over \lambda^2}=\lim_{k \rightarrow \infty}{ c_k \over c_{k+1}},
$$
where $\rho_k$ is the ratio of the smaller and larger eigenvalues of $D H_k(p_k)$. 

This universality can be explained rigorously if one shows that the {\it renormalization} operator
\begin{equation}\label{Ren}
R[F]=\Lambda^{-1}_F \circ F \circ F \circ \Lambda_F,
\end{equation}
where $\Lambda_F$ is some $F$-dependent coordinate transformation, has a fixed point, and the derivative of this operator is hyperbolic at this fixed point.

It has been argued in \cite{CEK2}  that $\Lambda_F$ is a diagonal linear transformation. Furthermore, such $\Lambda_F$ has been used in \cite{EKW1} and \cite{EKW2} in a computer assisted proof of existence of a reversible renormalization fixed point $F_*$ and hyperbolicity of the operator $R$. 

An exploration of a possible analytic  machinery has been undertaken in \cite{GK1} where it has been demonstrated that the fixed point $F_*$ is very close, in  some appropriate sense, to an area-preserving H\'enon-like map
\begin{equation}
H^*(x,u)=(\phi(x)-u,x-\phi(\phi(x)-u )),
\end{equation}
where $\phi$ solves the following one-dimensional problem of non-Feigenbaum type:
\begin{equation}\label{0_equation}
\phi(y)={2 \over \lambda} \phi(\phi(\lambda y)) -y.
\end{equation}

In this paper we will study the dynamics of the renormalization fixed point $F_*$  and maps in some neighbourhood of $F_*$. We would like to emphasize that the preservation of area for such maps leads to one important difference between the case at hand and the highly dissipative case of H\'enon-like maps: the former maps have no obvious invariant subsets in their domain. In particular, this means, that the construction of invariant Cantor sets carried out in \cite{dCLM}, which is using existence of invariant subsets in a crucial way, is not readily applicable to the case of area-preserving maps.

To construct hyperbolic sets, we will use the idea of covering relations (see, e.g. \cite{ZG04,Z09}) in rigorous computations. The Hausdorff dimension of the hyperbolic sets will be estimated with the help of the Duarte Distortion Theorem (see, e.g. \cite{Duarte1}) which enables one to use the distortion of a Cantor set to ultimately find bounds on the dimension. 

The structure of the paper is as follows: we begin by recalling the basic properties of area-preserving reversible maps in Section 2. In Section 3 we introduce some notation and recall some standard definitions from hyperbolic dynamics. In Section 4 we study the domain of analyticity of maps in a neighbourhood of the renormalization fixed point, and analytic continuation of the renormalization fixed point. Section 5 consists of the statements of our main theorems. In Section 6, we recall the definition and the main properties of the covering relations. In Section 7 we prove that any area-preserving reversible map in a neighbourhood of the renormalization fixed point has a transversal homoclinic orbit in its domain of analyticity. In Section 8 we construct a heteroclinic tangle for the renormalization fixed point. From its existence the existence of unbounded, and oscillating trajectories follow. In Section 9 we recall the Duarte Distortion Theorem. In Sections 10 and 11 we pr
 ove that the third iterate of any area-preserving reversible map in a neighbourhood of the renormalization fixed point has a horseshoe, and compute bounds on its Hausdorff dimension using the Duarte distortion theorem.

In a satellite paper \cite{GJ}, we prove that {\it infinitely renormalizable} maps in the neighbourhood of existence of the hyperbolic set for the third iterate also admit  a ``stable'' set. This set is a bounded invariant set, such that the maximal Lyapunov exponent for the third iterate is zero. In \cite{GJ} we provide an upper bound on the Hausdorff dimension of the stable set, and prove that the Hausdorff dimension is constant for all maps in some subset of infinitely renormalizable maps.

\section{Renormalization for area-preserving reversible maps} 

An ``area-preserving map'' will mean an exact symplectic diffeomorphism of a subset of ${\fR}^2$ onto its image.

Recall, that an area-preserving map can be uniquely specified by its generating function $S$:
\begin{equation}\label{gen_func}
\left( x \atop -S_1(x,y) \right) {{ \mbox{{\small \it  F}} \atop \mapsto} \atop \phantom{\mbox{\tiny .}}}  \left( y \atop S_2(x,y) \right), \quad S_i \equiv \partial_i S.
\end{equation}

Furthermore, we will assume that $F$ is reversible, that is 
\begin{equation}\label{reversible}
T \circ F \circ T=F^{-1}, \quad {\rm where} \quad T(x,u)=(x,-u).
\end{equation}

For such maps it follows from $(\ref{gen_func})$ that 
\begin{equation}\label{littles}
S_1(y,x)=S_2(x,y) \equiv s(x,y), 
\end{equation}

and
\begin{equation}\label{sdef}
\left({x  \atop  -s(y,x)} \right)  {{ \mbox{{\small \it  F}} \atop \mapsto} \atop \phantom{\mbox{\tiny .}}} \left({y \atop s(x,y) }\right).
\end{equation}

It is this ``little'' $s$ that will be referred to below as ``the generating function''. It follows from (\ref{littles}) that $s_1$ is symmetric. If the equation $-s(y,x)=u$ has a unique differentiable solution $y=y(x,u)$, then the derivative of such a map $F$ is given by the following formula:

\begin{equation}\label{Fder}
\hspace{-1.5cm}DF(x,u)=\left[ 
\begin{array}{c c}
-{s_2(y(x,u),x) \over s_1(y(x,u),x)} &  -{1 \over s_1(y(x,u),x)} \\
s_1(x,y(x,u))-s_2(x,y(x,u)) {s_2(y(x,u),x) \over s_1(y(x,u),x)}  & -{s_2(x,y(x,u)) \over s_1(y(x,u),x)} 
\end{array}
\right]. 
\end{equation}

We will now derive an equation for the generating function of the renormalized map $\Lambda_F \circ F \circ F \circ \Lambda_F^{-1}$.

Applying a reversible $F$ twice we get
$$
 \left({x'  \atop  -s(z',x')} \right) {{ \mbox{{\small \it  F}} \atop \mapsto} \atop \phantom{\mbox{\tiny .}}} \left({z' \atop s(x',z')} \right) = \left({z'  \atop  -s(y',z')} \right) {{ \mbox{{\small \it  F}} \atop \mapsto} \atop \phantom{\mbox{\tiny .}}} \left({y'\atop  s(z',y')} \right).
$$

It has been argued in \cite{CEK2}  that 

$$\Lambda_F(x,u)=(\lambda x, \mu u).$$

We therefore set  $(x',y')=(\lambda x,  \lambda y)$, $z'(\lambda x, \lambda y)= z(x,y)$ to obtain:

\begin{equation}\label{doubling}
\hspace{-1cm}\left(\!{x  \atop  -{ 1 \over \mu } s(z,\lambda x)} \!\right) \!{{ \mbox{{\small $\Lambda_F$}} \atop \mapsto} \atop \phantom{\mbox{\tiny .}}} \!\left(\!{\lambda x  \atop  -s(z,\lambda x)} \!\right) \!{{ \mbox{{\small \it  F $ \circ$ F}} \atop \mapsto} \atop \phantom{\mbox{\tiny .}}}\!\left(\!{\lambda y \atop s(z,\lambda y)}\! \right)   {{ \mbox{{\small \it  $\Lambda_F^{-1}$}} \atop \mapsto} \atop \phantom{\mbox{\tiny .}}} \left(\!{y \atop {1 \over \mu } s(z,\lambda y) }\!\right),
\end{equation}
where $z(x,y)$ solves
\begin{equation}\label{midpoint}
s(\lambda x, z(x,y))+s(\lambda y, z(x,y))=0.
\end{equation}

If the solution of $(\ref{midpoint})$ is unique, then $z(x,y)=z(y,x)$, and it follows from $(\ref{doubling})$ that the generating function of the renormalized $F$ is given by 
\begin{equation}
\tilde{s}(x,y)=\mu^{-1} s(z(x,y),\lambda y).
\end{equation}

One can fix a set of normalization conditions for $\tilde{s}$ and $z$ which serve to determine scalings $\lambda$ and $\mu$ as functions of $s$. For example, the normalization
$$s(1,0)=0$$ 
is reproduced for $\tilde{s}$ as long as 
$$z(1,0)=z(0,1)=1.$$

In particular, this implies that 
$$s(\lambda, 1)+s(0, 1)=0.$$

Furthermore, the condition
\begin{equation}\label{cond_s1}
\partial_1 s(1,0)=1
 \end{equation}
is reproduced as long as 
$$\mu=\partial_1 z (1,0).$$

We will now summarize the above discussion in the following definition of the renormalization operator acting on generating functions originally due  to the authors of \cite{EKW1} and \cite{EKW2}:

\medskip

\begin{definition}
\phantom{a}

\begin{eqnarray}\label{ren_eq}
\nonumber \\ {\cR}_{EKW}[s](x,y)=\mu^{-1} s(z(x,y),\lambda y),
\end{eqnarray}
where
\begin{eqnarray}
\label{midpoint_eq} 0&=&s(\lambda x, z(x,y))+s(\lambda y, z(x,y)), \\
0&=&s(\lambda,1)+s(0,1) \quad {\rm and} \quad \mu=\partial_1 z (1,0).
\end{eqnarray}
\end{definition}

\medskip

\begin{definition}
The Banach space of functions  $s(x,y)=\sum_{i,j=0}^{\infty}c_{i j} x^i y^j$, analytic on a bi-disk
$$ |x-0.5|<\rho, |y-0.5|<\rho,$$
for which the norm
$$\|s\|_\rho=\sum_{i,j=0}^{\infty}|c_{i j}|\rho^{i+j}$$
is finite, will be referred to as $\cA(\rho)$.
$\cA_s(\rho)$ will denote its symmetric subspace $\{s\in\cA(\rho) : s_1(x,y)=s_1(y,x)\}$.
\end{definition}

\medskip

As we have already mentioned, the following has been proved with the help of a computer in \cite{EKW1} and \cite{EKW2}:
\begin{theorem}\label{EKWTheorem}
There exist a polynomial $s_{\rm a} \in \cA_s(\rho)$ and  a ball $\cB_r(s_{\rm a}) \subset \cA_s(\rho)$, $r=6.0 \times 10^{-7}$, $\rho=1.6$, such that the operator ${\cR}_{EKW}$ is well-defined, analytic and compact on $\cB_r$. 

Furthermore, its derivative $D {\cR}_{EKW} \arrowvert_{\cB_r }$ has exactly two eigenvalues $\delta_1$ and $\delta_2$ of modulus larger than $1$, while 
$${\rm spec}(D {\cR}_{EKW} \arrowvert_{\cB_r }) \setminus \{\delta_1,\delta_2 \} \subset \{z \in \C: |z| \le \nu < 1\}.$$ 

Finally, there is an $s^* \in \cB_r$ such that
$$\cR_{EKW}[s^*]=s^*.$$
The scalings $\lambda_*$ and $\mu_*$ corresponding to the fixed point $s^*$ satisfy
\begin{eqnarray}
\label{lambda} \lambda_* \in [-0.24887681,-0.24887376], \\
\label{mu} \mu_* \in [0.061107811, 0.061112465].
\end{eqnarray}
 \end{theorem}

\medskip
 \begin{remark}
The radius of the contracting part of the spectrum ${\rm spec}(D {\cR}_{EKW}(s_*)  \setminus \{\delta_1,\delta_2 \}$ has been estimated in \cite{EKW2} to be $\nu=0.8$. 
 \end{remark}
\medskip

It follows from the above theorem that there exist a codimension $2$ local stable manifold $W^s_\loc(s^*)\subset \cB_r$.

\begin{definition}\label{Wdef}
A reversible map $F$ of the form (\ref{sdef}) such that $s\in W^s_\loc(s^*)$ is called infinitely renormalizable. The set of all reversible infinitely renormalizable maps is denoted by $\bW$.
\end{definition}

\section{Some notation and definitions}

We will use the following notation for the sup norm of a function $h$ and a transformation $H$ defined on 
some set $\cS \subset \fR^2 \! \quad {\rm or} \quad \! \fC^2$:
\begin{eqnarray}
\arrowvert h \arrowvert_\cS \equiv \sup_{(x,u) \in \cS }\{ |h| \},\\
\arrowvert H \arrowvert_\cS \equiv \max \{\sup_{(x,u) \in \cS }\{|\cP_x H| \}, \sup_{(x,u) \in \cS}\{|\cP_u H| \} \},
\end{eqnarray}
where $\cP_x$ and $\cP_u$ are projections on the corresponding components.

We will also use the notation $| \cdot |$ for the $l_2$ norm for vectors in $\fR^2$.


The interval enclosures of $\lambda_*$ and $\mu_*$ will be denoted 
\begin{eqnarray}
\label{blambda} \blambda &\equiv& [\lambda_-,\lambda_+]; \quad \lambda_-=-0.24887681, \quad \lambda_+=-0.24887376,\\
\label{bmu} \bmu & \equiv& [\mu_-,\mu_+]; \quad \mu_-=0.061107811, \quad \mu_+=0.061112465.
\end{eqnarray}

The corresponding interval enclosure for the linear map $\Lambda_*$ will be denoted $\bLambda$; if $(x,u) \in \fC^2$, then
\begin{equation}
 \bLambda(x,u) \equiv \left\{(\lambda x,\mu u) \in \fC^2: \lambda \in \blambda, \mu \in \bmu  \right\}. 
\end{equation}

The bound on the fixed point generating function $s^*$ will be called $\bs$:
\begin{equation}\label{bs}
\bs \equiv \left\{ s \in \cA_s(\rho) : \|s-s_{\rm a}\|_{\rho} \le  r=6.0 \times 10^{-7}  \right\},
\end{equation}
while the bound on the renormalization fixed point $F_*$ will be referred to as  $\bF$:
\begin{equation}\label{bF}
\bF \equiv \left\{ F: (x,-s(y,x)) \mapsto (y,s(x,y)) : s \in \bs \right\},
\end{equation}
where $s_{\rm a}$ is as in Theorem $\ref{EKWTheorem}$; the third iterate of this bound will be referred to as $\bG$. With $D\bG: p \mapsto [D\bG(p)]$ we denote an interval matrix valued function such that 
$$
[DG(p)]_{ij} \in [D\bG(p)]_{ij}, \quad \textrm{for all } G\in \bG,\,p\in\cD_3,
$$
where $\cD_3$ is the domain of $\bG$, and the bound on the operator norm of $DG$ for $G=F \circ F \circ F$, $F\in \bF$ on a set $\cS$ will be denoted
$$\| D\bG \|_\cS \equiv \sup_{F \in \bF }\left\{ \| D (F \circ F \circ F) \|_{\cS}\right\}.$$

Given a non-empty open set $\cD \subset \fC^n$ we will denote by $\cO_R(\cD)$ the set of reversible area-preserving maps $F: \cD \mapsto \fC^n$, analytic on $\cD$.

We will proceed with a collection of classical notations (see, eg, \cite{KH}) relevant to our following discussion.

\medskip


\medskip

\begin{definition} \label{hyperbolicset} (\underline{Hyperbolic set})
Let $\cM$ be a smooth manifold, and let $F$ be a diffeomorphism of an open subset $\cU \subset \cM$ onto its image.

A set $\cC$ is called hyperbolic for the map $F$ if there is a Riemannian metric on a neighbourhood $\cU$ of $\cC$, and $\beta<1<\delta$, such that for any $p \in \cC$ and $n \in \fN$ the tangent space $T_{F^n(p)}\cU$ admits a decomposition in two invariant subspaces:
$$T_{F^n(p)}\cM  =  E^+_n \oplus E^-_n, \quad 
DF(F^n(p)) E^\pm_n = E^\pm_{n+1},$$
on which the sequence of differentials is hyperbolic:
$$\|DF(F^n(p))\arrowvert_{E^-_n}\| < \beta, \quad \| DF^{-1} (F^n(p)) \arrowvert_{E^+_{n+1}} \| <  \delta^{-1}.$$
\end{definition}

\medskip

\begin{definition}\label{maximalset} (\underline{Locally maximal hyperbolic set})
Let $\cC$ be a hyperbolic set for $F: \cU \mapsto \cM$. If there is a neighbourhood $\cV$ of $\cC$ such that 
$\cC=\cap_{n \in \fZ} F^n(\bar{\cV}),$
then $\cC$ is called locally maximal.
\end{definition}

\medskip

\begin{definition} \label{shift}(\underline{Bernoulli shift})
Let  $\{0,1,\ldots,N-1\}^{\fZ}$ be the space of all two-sided sequences of $N$ symbols:
$$\{0,1,\ldots,N-1\}^{\fZ} = \{\omega=(\ldots,\omega_{-1},\omega_0,\omega_1,\ldots): \omega_i \in \{0,1,\ldots,N-1\}, i\in\fZ   \},$$
Define the Bernoulli shift on  $\{0,1,\ldots,N-1\}^{\fZ}$ as
$$\sigma_N (\omega)=\omega', \quad \omega'_n=\omega_{n+1}.$$
\end{definition}

\medskip

\begin{definition}\label{chain} (\underline{Topological Markov chain})
Let $A=(a_{i,j})_{i,j=0}^{N-1}$ be an $N\times N$ matrix whose entries are either $0$ or $1$. Let
 $$\{0,1,\ldots,N-1\}^{\fZ}_A=\{\omega \in  \{0,1,\ldots,N-1\}^{\fZ}: a_{\omega_n,\omega_{n+1}}=1, n \in \fZ \}.$$

The restriction 
$$\sigma_N \arrowvert_ {\{0,1,\ldots,N-1\}^{\fZ}_A} \equiv \sigma_A$$ 
is called a topological Markov chain determined by $A$.
\end{definition}

\medskip

\begin{definition}\label{homo_point}(\underline{Homoclinic and heteroclinic points})
Let $F : \cX \mapsto \cX$ be a homeomorphism on a metric space $(\cX,d)$. A point $x \in \cX$ is said to be {\it homoclinic} to the point $y \in \cX$ if 
$$\lim_{|n| \rightarrow \infty} d(F^n(x),F^n(y))=0.$$

A point $x$ is said to be {\it heteroclinic} to points $y_1$ and $y_2$ if 
$$\lim_{n \rightarrow \infty} d(F^n(x),F^n(y_1))=\lim_{n \rightarrow -\infty} d(F^n(x),F^n(y_2))=0.$$
In this case we will also say that there exists a {\it heteroclinic orbit} between points $y_1$ and $y_2$.

If $\cM$ is a differentiable manifold, $F \in {\rm Diff}^1(\cM)$ and $x \in \cM$ is a hyperbolic fixed point of $F$, then we say that $q \in \cM$ is a {\it transversal homoclinic} point to $x$ if it is a point of transversal intersection of the stable and unstable manifolds of $x$.

\end{definition}

\medskip

\begin{definition} \label{HDim}
Let $\cX$ be a metric space. If $\cS \subset \cX$, and $d \in [0,\infty)$, the $d$-dimensional \underline{Hausdorff content} of $\cS$ is defined as 
\begin{equation}\label{Hcont}
C_d^H(\cS)={\rm inf} \left\{\sum_i {\rm diam}(\cS_i)^d: \{\cS_i\} \, {\rm is \, a \, cover \, of \, } \cS \right\}.
\end{equation}
The \underline{Hausdorff dimension} of $\cS$ is defined as 
\begin{equation}
{\rm dim}_H(\cS) = {\rm inf} \left\{d \ge 0: C_d^H(\cS)=0\right\}.
\end{equation}
\end{definition}

\medskip

\section{Domain of analyticity of $\bF$ and $\bG$}

$F$ is defined implicitly by the generating function $s$, its domain is 
given as in \cite{EKW2}:
\begin{equation}\label{sDom}
\cD_s=\{(x,y) \in \fC^2 \,: \, |x-0.5|<1.6, |y-0.5|<1.6\}.
\end{equation}
To find the domain of $F\in\bF$, we note that its second argument is equal to $-s(y,x)$, for some $s\in\bs$ (see $(\ref{sdef})$). 
Thus, the domain of $F$, $\cD$, is given by:
\begin{equation}\label{fDom}
\!\!\!\!\!\!\!\!\!\!\!\!\!\!\!\!\!\!\!\!\!\!\!\!\!\!\!\!\!\!\!\!\!\!\!\!\!\!\!\!\!\! \cD = \{(x,u) \in \fC^2 \,:\, u=-s(y,x), \, s \in  \bs(y,x), \, |x-0.5|<1.6,\, |y-0.5|<1.6 \}.
\end{equation}

We denote by
\begin{equation}\label{fDomRe}
\tilde{\cD} =\{(x,u) \in \cD \,: \Im{x}=\Im{u}=0\},
\end{equation}
the real slice of $\cD$.

To solve nonlinear equations on the computer we use the interval Newton operator, see e.g. \cite{Ne90}.
\begin{definition}\label{NOp}(\underline{Interval Newton Operator}) Let $F:\cD \subseteq \fR^n \rightarrow \fR^n$, and let $D{\bf F}$ be an interval matrix valued function such that $[DF(x)]_{ij} \in [D({\bf F})(x)]_{ij},$ for all $x\in \cD$. Let $\bX \subset \cD$ be a Cartesian product of finite intervals, $\hat{x}\in\bX$, and assume that if $A\in D{\bf F}(\bX)$, then $A$ is non-singular. We define the interval Newton operator as:
$$
N(F,\bX,\hat{x})=\hat{x}-(D{\bf F})^{-1}(\bX){\bf F}(\hat{x}).
$$ 
The main properties of $N$ is that if $N(F,\bX,\hat{x}) \subset {\rm int}\bX$, then there exists a unique solution to $F(x)=0$ in $\bX$, which is contained in $N(F,\bX,\hat{x})$, and if $N(F,\bX,\hat{x})\cap\bX=\emptyset$, then there is no solution to $F(x)=0$ in $\bX$. 
\end{definition}

\begin{lemma}\label{lfDomRe}
There exists a non-empty open set $\bar{\cD}$, 
$$\cP_x \bar{\cD} \subset \{x \in \fR: |x-0.5| <1.6\},$$
such that for every $(x,u) \in \bar{\cD}$ there exists a unique solution of the interval Newton operator of the function $h_{(x,u)}(y)=u+s(y,x)$, $s \in \bs(y,x)$, that satisfies $|y-0.5|<1.6$.
\end{lemma}
\noindent {\it Proof.}
Set $\bY=\{y\in \fR \,:\,|y-0.5|<1.6\}$. Given $(x,u)$, let $N(h_{(x,u)},\bY,\hat{y})$ be the interval Newton operator (for some appropriately chosen $\hat{y} \in \bY$). We have verified that there exists a non-empty set $\bar{\cD}$, such that for all $(x,u) \in \bar{\cD}$
$$
N(h_{(x,u)},\bY,\hat{y}) \subset \bY.
$$
This verification is implemented in the program {\tt findDomain} of \cite{JP} It follows, see e.g. \cite{Ne90}, that there is a unique $y$, such 
that $(x,y)$ is in the real slice of $\cD_s$ and $h_{(x,u)}(y)=0$. Thus, $y$ is defined as a function of $(x,u)$, and by (\ref{fDom}) $(x,u)\in\tilde \cD$.   
\newline
{\flushright {$\Box$}}

Clearly, $\bar{\cD} \subset \tilde{\cD}$.

The generating function $s$ is analytic on a bi-disk, preferably $F$ 
should have a similar property. From the above construction, we can at least show that $F$ is defined on a complex neighbourhood of $\bar{\cD}$ from the above Lemma \ref{lfDomRe}:

\begin{figure}[t]
\begin{center}
\includegraphics[width=0.7 \textwidth]{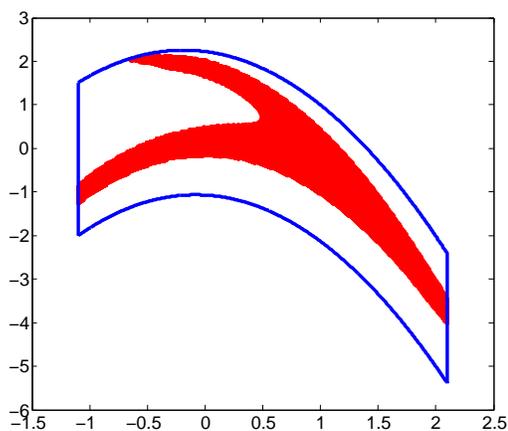}
\caption{The real slices of the domains of $F_*$ (blue) $G_*$ (red).}\label{DomPic}
\end{center}
 \end{figure}

\begin{lemma}
$\cD$ contains an open complex neighbourhood of the set $\bar{\cD}$. 
\end{lemma}
\noindent {\it Proof.}
$N(h_{(x,u)},y)$ is only well defined if $\bs_1(y,x) \neq 0$, so that the 
condition of the implicit function theorem is automatically satisfied on 
the solution set of $N(h_{(x,u)},\bY,\hat{y})$. Since 
$\bs_1$ is analytic on the bi-disk $\{|x-0.5|<1.6\} \times 
\{|y-0.5|<1.6\}$, there is an open neighbourhood in $\fC^2$ of the 
solution set where $\bs_1(y,x) \neq 0$. It follows from the implicit 
function theorem that $\bF$ is analytic on this neighbourhood.
\newline
{\flushright {$\Box$}}

Generally, we will denote the third iterate of a map $F \in \cO_R(\cD)$ as 
$G$. The domain of  $G$, $\cD_3$, is given by
\begin{equation}\label{gDom}
\cD_3=F^{-2} (\cD'' \cap \cD),
\end{equation}
where
$$\cD'' = F(\cD' \cap \cD), \quad \cD'= F(\cD).$$ 

\begin{figure}[t]
\begin{center}
\includegraphics[angle=-90,width=\textwidth]{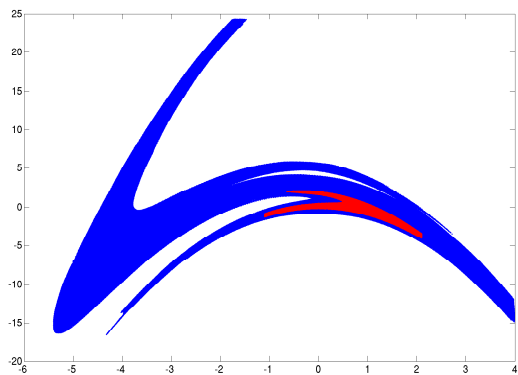}
\caption{The real slices of the domains $\cD^0_3$ (red) and  $\cD^1_3$ (blue).}\label{DomPicR}
\end{center}
 \end{figure}

Using the program {\tt findDomain}, we have verified that the real slice $\tilde{\cD}_3$ of $\cD_3$ is an open non-empty set. Approximations of $\tilde \cD$ and $\tilde{\cD_3}$ are shown in Figure \ref{DomPic}. We note that by using the renormalization equation $F_*=R[F_*]$, $F_*$ has, for any $k$, an analytic continuation to domains 
$$\cD^{k} \equiv \Lambda_*^{-1} (F_*^{-1}(F_*(\cD^{k-1})\cap  \cD^{k-1})), \quad \cD^0 \equiv \cD,$$ 
while $G_*$ has an analytic continuation to
$$\cD^{k}_3 \equiv \Lambda_*^{-1} (G_*^{-1}(G_*(\cD^{k-1}_3)\cap  \cD^{k-1}_3)), \quad \cD^0_3 \equiv \cD_3.$$ 
The real slices of the domains $\cD^0_3$ and $\cD^1_3$ are given in the Figure \ref{DomPicR}.

We will conclude this Section with the following

\begin{lemma}
All maps $ F \in \bF$ posses a hyperbolic fixed point $p_0=p_0(F)$, such that
\begin{itemize}
\item[1)] $\cP_x p_0 \in ( 0.57761843, 0.57761989)$, and $\cP_u p_0=0$;
\item[2)] $D F(p_0)$ has two real eigenvalues: 
\begin{eqnarray}
\nonumber e_+ &\in& (-2.05763559,-2.05759928),\\
\nonumber e_- &\in& (-0.48601715,-0.48598084).
\end{eqnarray} 
\end{itemize}
 \end{lemma}
\noindent {\it Proof.}
The bound on the fixed point has been obtained with the help of the interval Newton method. Hyperbolicity has been demonstrated by computing a bound on $D F$ using formula $(\ref{Fder})$.
\newline
{\flushright {$\Box$}}

\section{Statement of the main results}

We will now summarize our main findings.

Our first theorem describes the two-sided heteroclinic tangle for the fixed point map $F_*$. 

\begin{mainthm}\label{MTHM1}
The renormalization fixed point $F_*$ has the following properties:
\begin{itemize}
\item[1)] $F_*$ possesses a point $p_\pitchfork$ which is transversally homoclinic to the fixed point $p_0$; 
\item[2)] there exists a positive integer $n$ such that for any negative integer $k$ the map $F_*^n$ has a heteroclinic orbit $\cO_k$ between the periodic points  $\Lambda_*^{k}(p_0)$ and $\Lambda_*^{k+1}(p_0)$, and for any positive integer $k$ the map $F_*^{n \cdot 2^k}$ has a heteroclinic orbit $\cO_k$ between the periodic points  $\Lambda_*^{k}(p_0)$ and $\Lambda_*^{k-1}(p_0)$;
\item[3)]for any  $N \in  \fN$ and $M \in \fN$  there  exists an integer $n$ and an invariant set $\cZ_{F^n_*}$, such that 
$$F^n_* \arrowvert_{\cZ_{F^n_*}} \, { \approx \atop  \mbox{{\small \it  homeo}}   }  \, \sigma_A \arrowvert_{ \{-N,-N+1, \ldots,M-1,M\}^{\fZ}},$$
where $\sigma_A$ is the topological Markov chain defined by a $(N+M+1) \times (N+M+1)$  tridiagonal matrix
$$A=\left[ \
\begin{tabular}{ccccccccc} 
1 & 1 &0 & 0 &\ldots & 0& 0 & 0& 0 \\
1 & 1 &1 & 0 &\ldots & 0& 0 & 0& 0 \\
0 & 1 &1 & 1 &\ldots & 0& 0 & 0& 0 \\
. & . &. & . &\ldots & .& . & .& .\\
0 & 0 &0 & 0 &\ldots &1 & 1 & 1 & 0 \\
0 & 0 &0 & 0 &\ldots &0 & 1 & 1& 1 \\
0 & 0 &0 & 0 &\ldots &0 & 0 & 1& 1
 \end{tabular}  \right];
$$
\item[4)] for any $R>\varepsilon>0$ there exists a point $p \in \cD$ and an $n \in \fN$  such that $|F^n_*(p)|>R$ (\underline{an unbounded orbit}), and a point $q \in \cD$ and $n,m \in \fN$  such that $|F^n_*(q)|>R$ and $|F^m_*(q)|<\varepsilon$ (\underline{an oscillating orbit}).
\end{itemize}
\end{mainthm}

Our second result demonstrates that all locally infinitely renormalizable maps admit a hyperbolic set in their domain of analyticity.

\begin{mainthm}\label{MTHM2}
Any $F \in \bF$ admits a hyperbolic set $\cC_G \subset \cD_3$ for $G \equiv F \circ F \circ F$; 
$$G \arrowvert_{\cC_G} \, { \approx \atop  \mbox{{\small \it  homeo}}   }  \, \sigma_2 \arrowvert_{\{0,1\}^{\fZ}},$$ 
whose Hausdorff dimension satisfies:
$$0.7673 \ge  {\rm dim}_H(\cC_G) \ge \varepsilon,$$
where $\varepsilon \approx 0.00013 \, e^{-7499}$ is strictly positive.  
\end{mainthm}

\section{Topological tools}
The main tools of our proofs are {\it covering relations} 
\cite{Z97,ZG04} and {\it cone conditions} \cite{KWZ07,Z09}, see also 
\cite{GZ01}
for proving the existence of homoclinic and heteroclinic orbits. To make 
the present paper reasonably self-contained we include a brief 
introduction to the necessary concepts.

\subsection{H-sets and Covering relations}
The notion of an \textit{h-set} and a covering relation first appeared 
in \cite{Z97}, the most thorough treatment is \cite{ZG04}. The basic 
idea is to construct computable conditions for the existence of a 
semi-conjugacy to symbolic dynamics. This is done by constructing 
h-sets, i.e. hyperbolic-like sets, that cross each other in a 
(topologically) non-trivial way. We denote by $B_n(c,r)$, the open ball 
in $\mathbb{R}^n$ with centre $c$ and radius $r$, and $S^n(c,r)=\partial 
B_{n+1}(c,r)$.

\begin{definition}\label{hset} 
An \textit{h-set} is a quadruple consisting of
\begin{itemize}
 \item a compact subset $|N|$ of $\mathbb{R}^n$,
\item a pair of numbers $u(N), s(N)\in \{0,1,2,\dots \},$ with 
$u(N)+s(N)=n$,
\item a homeomorphism $c_N: 
\mathbb{R}^n\rightarrow\mathbb{R}^n=\mathbb{R}^{u(N)}\times 
\mathbb{R}^{s(N)}$, such that
$$c_N(|N|)=\overline{B_u(N)}(0,1)\times\overline{B_s(N)}(0,1).$$
\end{itemize}
\end{definition}
 
We denote such a quadruple by $N$. We usually drop the bars on the 
support and refer to it as $N$. Furthermore,
$$N_c=\overline{B_{u(N)}}(0,1)\times\overline{B_{s(N)}}(0,1),$$
$$N_c^-=\partial\overline{B_{u(N)}}(0,1)\times\overline{B_{s(N)}}(0,1),$$
$$N_c^+=\overline{B_{u(N)}}(0,1)\times\partial\overline{B_{s(N)}}(0,1),$$
$$N^-=c_N^{-1}(N_c^-), \quad N^+=c_N^{-1}(N_c^+).$$

$u(N)$ and $s(N)$ are the nominally unstable and stable directions, 
respectively. The idea of a covering relation between two h-sets is that 
the image of the first should be mapped transversally across $N^-$ ({\it the exit set of $N$}) and inside of $N^+$ ({\it the entrance set of $N$}).  We formally define it.

\begin{definition}
Assume $N,M$ are h-sets, such that $u(N)=u(M)=u$ and $s(N)=s(M)=s$. Let 
$f:N\rightarrow \mathbb{R}^n$ be a continuous map. Let $f_c=c_M\circ 
f\circ c_N^{-1}:N_c\rightarrow \mathbb{R}^u\times \mathbb{R}^s$. We say 
that $$N\cover{f}M$$
($N$ $f-$covers $M$) iff the following conditions are satisfied

1. There exists a continuous homotopy $h: [0,1]\times N_c \rightarrow 
\mathbb{R}^u \times \mathbb{R}^s$, such that the following conditions 
hold true
$$ h_0=f_c, $$
$$h([0,1],N_c^-)\cap M_c =\emptyset, $$
$$h([0,1],N_c)\cap M_c^+=\emptyset.$$

2. There exists a linear map $A:\mathbb{R}^u \rightarrow \mathbb{R}^u$, 
such that
$$h_1(p,q)=(Ap,0) \quad \textrm{ for } p\in \overline{B_u}(0,1)\quad 
\textrm{ and } q\in \overline{B_s}(0,1),$$
$$A(\partial B_u(0,1))\subset \mathbb{R}^u\setminus \overline{B_u}(0,1).$$

$h_1$ is called a model map for the relation $N\cover{f}M$.
\end{definition}

The maps that we study in the present paper are reversible, which we use 
to reduce the amount of computation. For such maps the following 
definition is useful, see also the discussion in \cite{ZG04}.

\begin{definition}
Let $N$ be an $h$-set. We define the $h$-set $N^T$ as follows:
\begin{itemize}
 \item The compact subset of the quadruple $N^T$ is the compact subset 
of the quadruple $N$, also denoted by $N$,
\item $u(N^T)=s(N),$ $s(N^T)=u(N)$.
\item The homeomorphism $c_{N^T}: \mathbb{R}^n\rightarrow 
\mathbb{R}^n=\mathbb{R}^{u(N^T)}\times \mathbb{R}^{s(N^T)}$ is defined by
$$c_{N^T}(x)=j(c_N(x)),$$
where $j:\mathbb{R}^{u(N)}\times \mathbb{R}^{s(N)} \rightarrow 
\mathbb{R}^{s(N)}\times \mathbb{R}^{u(N)}$ is given by $j(p,q)=(q,p)$.
\end{itemize}
\end{definition}
 
\begin{definition}
Assume N,M are $h$-sets, such that $u(N)=u(M)=u$ and $s(N)=s(M)=s$. Let 
$g:\Omega\subset\mathbb{R}^n \rightarrow \mathbb{R}^n$. Assume that 
$g^{-1}:M\rightarrow \mathbb{R}^n$ is well defined and continuous. We 
say that $N\bcover{g}M$ ($N$ $g-$backcovers $M$) iff $M^T\cover{g^{-1}}N^T$.
\end{definition}

\begin{definition}
Let $N$ and $M$ be $h$-sets. We say that $N$ generically $f$-covers $M$ 
($N\gcover{f}M$) if $N\cover{f}M$ 
or $N\bcover{f}M$.
\end{definition}

The main property of covering relations is contained in the following 
theorem \cite[Corollary 7]{ZG04}.
\begin{theorem}\label{tCovRel}
Assume that we have the following chain of covering relations:
$$N_0 \gcover{f_1}N_1 \gcover{f_2} N_2 
\Longleftrightarrow \dots \gcover{f_k} N_k,$$
then there exists a point $x\in {\rm int }\, N_0$, such that
$$f_i\circ f_{i-1}\circ \dots \circ f_1(x)\in {\rm int }\, N_i, \quad 
i=1,\dots,k.$$
Moreover, if $N_k=N_0$, then $x$ can be chosen so that
$$f_k\circ f_{k-1}\circ\dots\circ f_1(x)=x.$$
\end{theorem}
In our proofs the hypothesis of Theorem \ref{tCovRel} is verified using the routine {\tt checkCoveringRelations} of \cite{JP}.

\subsection{Cone conditions for h-sets}
Theorem \ref{tCovRel} gives a computational tool to prove the existence 
of orbits with prescribed symbolic dynamics. To prove that such orbits are unique 
one would ideally require hyperbolicity of the map in a neighbourhood of the orbit.
Typically, this is proved by constructing invariant cone fields. An alternative method to prove uniqueness is provided by  
covering relations with cone conditions, first described in \cite{KWZ07}, the method is studied in further detail in 
\cite{Z09}, which we follow below.

 \begin{definition}
Let $N\subset \mathbb{R}^n$ be an $h$-set and 
$Q:\mathbb{R}^n\rightarrow\mathbb{R}^n$  be a quadratic form
\begin{equation}
Q((x,y))=\alpha(x)-\beta(y),\quad(x,y)\in\mathbb{R}^u\times\mathbb{R}^s,
\end{equation}
where $\alpha:\mathbb{R}^{u(N)}\rightarrow\mathbb{R}$, and 
$\beta:\mathbb{R}^{s(N)}\rightarrow\mathbb{R}$ are positive definite 
quadratic forms.
 \end{definition}

The pair $(N,Q)$ is called an \textit{$h$-set with cones}.

\begin{definition}
Assume that $(N,Q_N)$ and $(M,Q_M)$ are $h$-sets with cones, such that 
$u(N)=u(M)=u$, and let $f:N\rightarrow \mathbb{R}^{dim(M)}$ be 
continuous. Assume that $N\cover{f}M$. We say that f 
satisfies the cone condition (with respect to the pair $(N,M)$) iff for 
any $p_1, p_2\in N_c$, $p_1\neq p_2$ holds
\begin{equation}
Q_M(f_c(p_1)-f_c(p_2)) > Q_N(p_1-p_2).
\end{equation}
\end{definition}

\begin{definition}
Assume that $(N,Q_N)$ and $(M,Q_M)$ are $h$-sets with cones, such that 
$u(N)=u(M)=u$ and $s(N)=s(M)=s$, and let $f:N\rightarrow 
\mathbb{R}^{u+s}$ be continuous. Assume that $N\bcover{f}M$. 
We say that f satisfies the cone condition (with respect to the 
pair $(N,M)$) iff for any $q_1, q_2\in M_c$, $q_1\neq q_2$ holds
\begin{equation}
Q_M(q_1-q_2) > Q_N(f_c^{-1}(q_1)-f_c^{-1}(q_2)).
\end{equation}
\end{definition}


\begin{remark}\label{rSMfd}
Note, for details see the discussion in \cite{Z09}, that stable 
manifolds are propagated backwards through a chain of covering relations 
with cone conditions. Similarly, unstable manifolds are propagated forwards through a chain of covering relations.
In the planar case, used in the present paper, the 
Lipschitz constant of the stable manifold is given by 
$\sqrt{\frac{\beta}{\alpha}}$, where $alpha$ and $beta$ are the constant
coefficients of the corresponding quadratic forms.  
\end{remark}

To verify that the cone conditions hold, we use the following Lemma 
\cite[Lemma 8]{Z09}.
\begin{lemma}\label{lQF}
Assume that for any $B\in D{\bf f}_c(N_C)$, the quadratic form
\begin{equation}
V(x)=Q_M(Bx)-Q_N(x)
\end{equation}
is positive definite, then for any $p_1,p_2\in N_c$ such that $p_1\neq p_2$
\begin{equation}
Q_M(f_c(p_1)-f_c(p_2))>Q_N(p_1-p_2).
\end{equation}
\end{lemma}
In our proofs the hypothesis of Lemma \ref{lQF} is verified using the routine {\tt checkConeConditions} of \cite{JP}.

\begin{remark}\label{rSep}
The proof of Lemma \ref{lQF} also demonstrates that 
$$Q_M(f_c(p_1)-f_c(p_2))-Q_N(p_1-p_2)=O(|p_1-p_2|^2),$$
and yields a uniform lower bound on $O(|p_1-p_2|^2)$ for each pair of $h$-sets:
$$Q_M(f_c(p_1)-f_c(p_2))-Q_N(p_1-p_2)\ge \epsilon_{MN}|p_1-p_2|^2.$$
\end{remark}

\begin{theorem}\label{tCovRelCones}
Let $A=(a_{ij})_{i,j=0}^{k-1}$ be a $k\times k$ matrix whose entries are either $0$ or $1$. 
Assume that the following set of covering relations with cone conditions holds 
$$N_i\gcover{f}N_j \quad {\rm if }\, a_{ij}=1.$$
Then, each sequence in $\{0,1,\ldots,k-1\}^{\fZ}_A$ (see Definition \ref{chain}) is realised by a unique orbit.
\end{theorem}
\noindent {\it Proof.}
Assume that $p_1 \in N_{\omega_0}$ and $p_2 \in N_{\omega_0}$, are such that $p_1 \neq p_2$, $f^i(p_1)\in N_{\omega_i}$ and $f^i(p_2)\in N_{\omega_i}$, for all $i \in\fZ$, where $\omega = (...,\omega_{-1},\omega_0, \omega_{1},...)\in  \{0,1,\ldots,k-1\}^{\fZ}_A$. Let $\epsilon=\min_{\{(i,j): a_{ij}=1\}}{\epsilon_{N_{j}N_{i}}}$, where ${\epsilon_{N_{j}N_{i}}}$ is as in Remark \ref{rSep}. There are three cases. 

(i) Let $Q_{N_{\omega_0}}(p_1-p_2)>0$, then we have that for any $i>0$
\begin{eqnarray}
\nonumber \hspace{-2.4cm}Q_{N_{\omega_i}}(f_c^i(p_1)-f_c^i(p_2))&>&Q_{N_{\omega_{i-1}}}(f_c^{i-1}(p_1)-f_c^{i-1}(p_2))+\epsilon |f_c^{i-1}(p_1)-f_c^{i-1}(p_2)|^2 \\
\nonumber &>& \cdots \\
& > & Q_{N_{\omega_{0}}}(p_1-p_2) +\epsilon \sum_{k=0}^{i-1} |f_c^{k}(p_1)-f_c^{k}(p_2)|^2.
\end{eqnarray}
Since by assumption $f^i(p_1)-f^i(p_2)$ stays inside a compact neighbourhood of the origin, for all $i\in\fZ$, where $Q_{N_{\omega_i}}$ is bounded from above, the left hand side is bounded. Hence $f_c^i(p_1)-f_c^i(p_2) \rightarrow 0$, and therefore $Q_{N_{\omega_i}}(f_c^i(p_1)-f_c^i(p_2))\rightarrow 0$, but 
$Q_{N_{\omega_i}}(f_c^i(p_1)-f_c^i(p_2))> Q_{N_{\omega_0}}(p_1-p_2)>0$, since $p_1\neq p_2$. This is a contradiction. 

(ii) Let $Q_{N_{\omega_0}}(p_1-p_2)=0$, then $Q_{N_{\omega_1}}(f_c(p_1)-f_c(p_2))>0$, and the above argument holds for the pair $f_c(p_1)$ and $f_c(p_2)$.

(iii) Let $Q_{N_{\omega_0}}(p_1-p_2)<0$
we get the same contradiction for $i<0$ and ``$>$'' replaced by ``$<$''.
\newline
{\flushright {$\Box$}}

\section{A transversal homoclinic orbit for $\bF$}
The proofs of the theorems in this and the next sections were implemented in a 
computer program \cite{JP} using the software package CAPD \cite{CAPD} for 
interval arithmetic and covering relations. All computations were 
performed on an Intel Xeon 2.3GHz, 64bit processor with 16Gb of RAM. The 
program was compiled with \texttt{gcc}, version 4.1.2.

\begin{figure}[t]
\begin{center}
\includegraphics[width=0.7\textwidth]{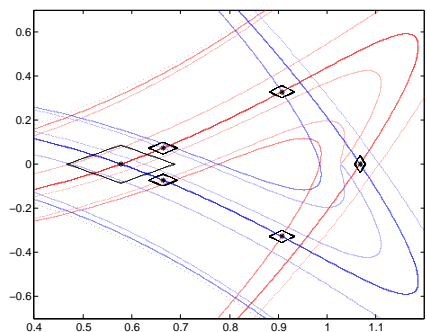}
\caption{The homoclinic tangle at the fixed point $p_0$. The stable manifold is given in blue, unstable --- in red. The black boxes correspond to the $h$-sets used in the proof.}\label{fHomoTraj}
\end{center}
\end{figure}

Numerical experiments indicate that for all $F\in\bF$ there exists a transversal 
homoclinic point $p_\pitchfork$, approximately located at 
$(1.067707,0)$, to the fixed point $p_0$, approximately located at 
$(0.577619,0)$, see Figure $\ref{fHomoTraj}$. We construct a sequence of covering relations and 
quadratic forms, and verify that the hypotheses of Theorem \ref{tCovRel} and Lemma 
\ref{lQF} are satisfied. Here, and in the rest of the paper, we use formula $(\ref{Fder})$ to provide bounds on the derivative of $F$.
\begin{theorem}\label{tHorse}
Any map $F\in\bF$ possesses a point $p_\pitchfork$  transversally homoclinic to $p_0$.
The map $F^{12}$ admits a locally maximal invariant set $\cC_{F^{12}} \ni p_*, p_\pitchfork$:
$$F^{12} \arrowvert_{\cC_{F^{12}}} { \approx \atop  \mbox{{\small \it  homeo}}   } \sigma_2 \arrowvert_{\{0,1\}^{\fZ} }.$$
\end{theorem}
\noindent {\it Proof.}
To avoid the ``flip'' of the system we study the second iterate of $F$, 
following the stable manifold from the homoclinic point to the fixed 
point. We construct four $h$-sets, $B_i$, $i=1,2,3,4$, as in Table \ref{tHomoTraj}: box $B_4$ around $p_0$,  box $B_1$ around $p_\pitchfork$ and two intermediate boxes $B_2$, and $B_3$.
$B_4$ is spanned by the unstable and stable eigenvectors of $DF(p_0)$. The other boxes are 
experimentally rotated until two sides are roughly parallel with the 
stable manifold. We prove (using the routines {\tt checkCoveringRelations} and {\tt  checkConeConditions} in \cite{JP}) that the following covering sequence with cone conditions holds 
$$B_1\Rightarrow B_2 \Rightarrow B_3 \Rightarrow B_4 \Rightarrow B_4.$$

The map $F$ is reversible, by reversing time and reflecting the $h$-sets 
in $\{u=0\}$ our cover is also a back-cover from $p_\pitchfork$ to 
$p_0$ along the unstable manifold of $p_0$. Therefore we also have the following covering relations with cone conditions

$$B_4\Leftarrow T(B_3) \Leftarrow T(B_2) \Leftarrow B_1.$$
Existence of a homoclinic orbit follows from Theorem \ref{tCovRelCones}. The transversality of the intersection follows from the fact that we use cones with Lipschitz constant less 
than one,  ${\rm Lip}=0.95$. 
\newline
{\flushright {$\Box$}}

Note that both $p_0$ and $p_\pitchfork$ depend analytically on $F$.
Part $1$ of the first Main Theorem follows immediately, since $F_*\in\bF$.
\begin{table}[h]
\begin{center}
\begin{footnotesize}
\begin{tabular}{c|cc}
Box \# & Centre & ``Unstable'' Direction, ${\bf e}^u$ \\ \hline
$B_1$ & $(1.067744, 0.0)$ & $(0.271460881389316,0.962449473933741)$ \\
$B_2$ & $(0.90733000,  -0.328353)$ & $(0.653774812595729,  0.691649449891044)$ \\
$B_3$ & $(0.66401000,  -0.0727302)$ & $(-0.736700425093906,  -0.676219257095255)$ \\
$B_4$ & $(0.5776305,  0.0)$ & $(0.788578889012330, 0.614933602760558)$
\end{tabular}
{\caption{Boxes $B_i$ used to prove the existence of a transversal 
homoclinic orbit. The size of each side of $B_4$ is $0.035|{\bf e}_4|$; the size of $B_i$, $i=1,2,3$, is $0.01 |{\bf e}^u|$. The ``stable'' spanning vector is defined as ${\bf e}^s=T({\bf e}^u)$. The covering sequence is  $B_1\Rightarrow B_2 \Rightarrow B_3 \Rightarrow B_4 \Rightarrow B_4$.}\label{tHomoTraj}}
\end{footnotesize}
\end{center}
\end{table}

\begin{remark}
The numbers describing boxes $B_4$ in the Table $\ref{tHomoTraj}$ are non-representable on the computer. We have used these numbers as the input for the routine {\tt  checkConeConditions}. During its execution, the numbers are rounded up to the nearest representable number; therefore, strictly speaking, Theorem  $\ref{tHorse}$ has been proved for boxes different from (although very close to) those reported in the table. Since the specific values of the true parameters are largely irrelevant, we will not go into pains of computing and reporting bounds on these parameters. This will be the situation through out the paper: all specific numerical data appearing in various tables and statements of theorems are (very accurate) approximations of the true data for which the results have been proved on the computer. 
\end{remark}


\section{A bi-infinite heteroclinic tangle for $F_*$}

In this Section we will demonstrate that the fixed point $F_*$ has a sequence of hyperbolic periodic points $\{z_k \}$, $k \in \fZ$,
 $$|z_k|  \converge{k \rightarrow \infty} 0  \quad {\rm and} \quad |z_k| \converge {k \rightarrow -\infty} \infty,$$
such that their stable and unstable manifolds form a heteroclinic tangle.

We first demonstrate that there exists a heteroclinic orbit between the fixed point $z_0 \equiv p_0$ and the period two point $z_1 \equiv \Lambda_*(p_0)$.  

We proceed as in the previous section:  we construct a covering sequence of $h$-sets, as in Table \ref{tP2HeteroTraj}, consisting of: one box $B_0$ around $p_0$, one box $B_1$ around $\Lambda_*(p_0)$ and $11$ intermediate boxes on the orbit. In this sequence, $B_0$ is defined as $\Lambda_*^{-1}(B_{1})$.

\begin{table}[h]
\begin{center}
\begin{footnotesize}
\begin{tabular}{c|ccccc}
$\!\!\!\!$ Box $\# \!\!\!$& Centre & $l^s_i$ & $l^u_i$  & ${\bf e}^s_i$ &  ${\bf e}^u_i$ \\ \hline
$\!\!\!B_1\!\!\!$ &  $\!\!\!(-0.143816, 0)\!\!\!$ & $\!\!\!0.00072\!\!\!$ & $\!\!\!0.00072\!\!\!$ &$\!\!\!(0.982114, 0.188285)\!\!\!$&   $\!\!\!(0.982114, -0.188285)\!\!\!$ \\
$\!\!\!B_2\!\!\!$ & $\!\!\!(-0.141389, -0.000460854)\!\!\!$ & $\!\!\!0.00054\!\!\!$ & $\!\!\!0.00015\!\!\!$ & $\!\!\!(0.982114, 0.188285)\!\!\!$ &  $\!\!\!(0.982114, -0.188285)\!\!\!$  \\
$\!\!\!B_3\!\!\!$ & $\!\!\!( -0.133448, -0.00190821)$ &$\!\!\!0.00045\!\!\!$&$\!\!\!0.00015\!\!\!$&$\!\!\!(0.982114, 0.188285)\!\!\!$&$\!\!\!(0.982114, -0.188285)\!\!\!$  \\
$\!\!\!B_4\!\!\!$ & $\!\!\!(-0.098565, -0.00713661)$ &$\!\!\!0.00047\!\!\!$&$\!\!\!0.0002\!\!\!$&$\!\!\!(0.987, 0.16072)\!\!\!$&$\!\!\!(0.997, -0.0774016)\!\!\!$ \\
$\!\!\!B_5\!\!\!$ & $\!\!\!(0.0582909, -0.00620757)\!\!\!$ &$\!\!\!0.00055\!\!\!$&$\!\!\!0.00015\!\!\!$&$\!\!\!(0.994, -0.10938)\!\!\!$&$\!\!\!(0.995, 0.0998749)\!\!\!$  \\
$\!\!\!B_6\!\!\!$ & $\!\!\!(0.106679, 0.0110872)\!\!\!$ &$\!\!\!0.0004\!\!\!$&$\!\!\!0.0001\!\!\!$&$\!\!\!(0.975, -0.222205)\!\!\!$&$ \!\!\!(0.982114, 0.188285)\!\!\!$  \\
$\!\!\!B_7\!\!\!$ & $\!\!\!(-0.2842, 0.070894)\!\!\!$ &$\!\!\!0.0004\!\!\!$&$\!\!\!0.0001\!\!\!$&$\!\!\!(0.84, 0.542586)\!\!\!$&$\!\!\!(0.866025, -0.5)\!\!\!$  \\
$\!\!\!B_8\!\!\!$ & $\!\!\!(0.533255, 0.0331362)\!\!\!$ &$\!\!\!0.0003\!\!\!$&$\!\!\!0.0002\!\!\!$&$\!\!\!(0.8, -0.6)\!\!\!$&$\!\!\!(0.788579, 0.614934)\!\!\!$  \\
$\!\!\!B_9\!\!\!$ & $\!\!\!(0.575172, 0.00190701)$ &$\!\!\!0.0003\!\!\!$&$\!\!\!0.0002\!\!\!$&$\!\!\!(0.788579, -0.614934)\!\!\!$&$\!\!\!(0.788579, 0.614934)\!\!\!$  \\
$\!\!\!B_{10}\!\!\!$ & $\!\!\!(0.577516, 0.000132348)\!\!\!$ &$\!\!\!0.0005\!\!\!$&$\!\!\!0.0002\!\!\!$&$\!\!\!(0.788579, -0.614934)\!\!\!$&$\!\!\!(0.788579, 0.614934)\!\!\!$  \\
$\!\!\!B_{11}\!\!\!$ & $\!\!\!(0.577619,0)\!\!\!$ & $\!\!\!0.001\!\!\!$ & $\!\!\!0.001\!\!\!$  & $\!\!\!(0.788579, -0.614934)\!\!\!$& $\!\!\!(0.788579, 0.614934)\!\!\!$ \\
\end{tabular}
\caption{The data used to prove the existence of a heteroclinic orbit 
for $F_*^2$ between $z_0$ and $z_1$. Vectors ${\bf e}^{u,s}_i$ are the ``stable'' and ``unstable'' spanning vectors of the rectangles. The length of the sides of the rectangles $B_i$ is $2 \cdot l^{u,s}_i \cdot |{\bf e}^{u,s}_i|$.  The covering sequence is
$B_1\Rightarrow B_1\Rightarrow B_2\Rightarrow B_3\Rightarrow B_4\Rightarrow B_5 
\Rightarrow B_6 \Rightarrow B_7 \Rightarrow B_8 \Rightarrow B_9\Rightarrow B_{10} \Rightarrow B_{11} \Rightarrow B_0=\Lambda_*^{-1}(B_1) \Rightarrow B_0=\Lambda_*^{-1}(B_1)$. The first and the last two covering relations are for $F^2_*$, the other ones --- for $F^4_*$.}\label{tP2HeteroTraj}
\end{footnotesize}
\end{center}
\end{table}

Since the orbit is on the tangle between the unstable and stable 
manifolds, there are parts of it where it is hard to construct expanding 
and contracting directions directly for the second iterate of the map. 
The fourth iterate, however, exhibits nice hyperbolic  like behaviour along a pseudoorbit. We use the programs {\tt checkCoveringRelations} and {\tt checkConeConditions} to verify that the hypotheses of Theorem  \ref{tCovRel} and Lemma \ref{lQF} are satisfied for  $F^4_*$ on the sequence of $h$-sets. In addition, we prove that they are satisfied for $F^2_*$ for the covering $B_1 \Rightarrow B_1$ (which implies that the covering relations and cone conditions are also satisfied for $B_0 \Rightarrow B_0$). We use cones with Lipschitz constant $0.9$ in the proof. 

The reversibility of $F_*$ implies that the existence of a heteroclinic orbit from $z_1$ to $z_0$ entails the existence of a 
heteroclinic orbit from $z_0$ to $z_1$. It also implies existence of homoclinic orbits for both points, that enters a small neighbourhood of 
the other one. The sizes of these neighbourhoods are given by the sizes 
of the boxes $B_0$ and $B_1$ in Table \ref{tP2HeteroTraj}. 

To construct the pseudoorbit we non-uniformly distribute points in the unstable eigendirection of $z_1$,
and locate an initial point that passes as close as possible to $z_1$ and $z_0$. 
The finest discretization of the approximate unstable vector for the period two point is done on the scale $10^{-14}$. 
The first point on the pseudoorbit is in the unstable eigendirection of $z_1$ 
at a  distance of $0.00058120854283$, and the last point is at a distance of $0.000167835$ from $z_0$.

We will now demonstrate the following eight covering relations:
$$ B_1 \cover{} B_0, \quad  B_1 \cover{} B_1; \quad B_0 \bcover{} B_1, \quad B_1 \bcover{} B_1;$$
$$B_0 \cover{} B_0, \quad B_0 \bcover{} B_0; \quad B_0 \cover{} B_1, \quad B_1 \bcover{} B_0.$$ 

\noindent {\it Step 1.} Set 
\begin{eqnarray}
\nonumber \hspace{-2.0cm}\phantom{a}_1B_1^0&\hspace{-1.0cm}=\!&\hspace{-0.2cm}\!F_*^{-4}(F^4_*(B_1) \cap F^{-4}_*(\\ & & \hspace{-0.2cm}\ldots \cap F^{-4}_*(F^4_*(B_8) \cap F^{-4}_*(F^4_*(B_9) \cap F^{-4}_*(F^4_*(B_{10}) \cap B_{11})))\ldots)),\\
\hspace{-2.0cm}\nonumber \phantom{a}_1B_1^1&\hspace{-1.0cm}=\!&\hspace{-0.2cm}\!F_*^{-4}(F^4_*(B_1) \cap F^{-4}_*( \\ & & \hspace{-0.2cm}\ldots \cap F^{-4}_*(F^4_*(B_1) \cap F^{-4}_*(F^4_*(B_1) \cap F_*^{-4}(F^4_*(B_{1}) \cap B_{1})))\ldots))
\end{eqnarray} 

Clearly, the iterate $\hF_* \equiv F^{42}_*$ is well defined and analytic on $\Bm{1}{0} \cup \Bm{1}{1}$, and 
\begin{equation}\label{cov_rel1}
\Bm{1}{0} \cover{\hF_*} B_0, \quad  
\Bm{1}{1} \cover{\hF_*} B_1.
\end{equation} 

The set $B_1$ is a parallelogram. Let $\pi^s_1$ be a projection on the span of ${\bf e}^s_1$ ``along'' the vector ${\bf e}^u_1$, and  $\pi^u_1$ be a projection on the span of ${\bf e}^u_1$ ``along'' the vector ${\bf e}^s_1$. Then the sets $\Bm{1}{i}$, $i=0,1$, are ``vertical'' in the sense that
 $\pi^s_1(\Bm{1}{i})=\pi^s_1(B_1)$, while the sets $\hF_*(\Bm{1}{i})$, $i=0,1$, are ``horizontal'' in the sense that
  $\pi^u_0\hF_*(\Bm{1}{0})=\pi^u_0(B_0)$ and $\pi^u_1\hF_*(\Bm{1}{1})=\pi^u_1(B_1)$.

\medskip

\noindent {\it Step 2.} Since the set $B_1$ is symmetric: $T(B_1)=B_1$, we have that 
\begin{equation}\label{cov_rel2}
T(\Bm{1}{0}) \cover{\hF^{-1}_*} B_0, \quad T(\Bm{1}{1}) \cover{\hF^{-1}_*} B_1.
\end{equation}

\medskip
 
\noindent {\it Step 3.} We get from the period doubling equation:
$$\hF_*(\Lambda^{-1}_*(\Bm{1}{1})) =\Lambda^{-1}_*(\hF_*^2(\Bm{1}{1})).$$ 
 Clearly, since $\hF_*$ is well-defined and analytic on $\Bm{1}{1}$, and $\Bm{1}{1} \cover{\hF_*} \Bm{1}{1}$  there exists a vertical subset $\Bm{2}{1} \subset \Bm{1}{1}$ on which $\hF^2_*$ is analytic. Set $\BmZ{1}{0}=\Lambda^{-1}_*(\Bm{2}{1})$. Then, since $B_0=\Lambda^{-1}_*(B_1)$, $\BmZ{1}{0}$ covers $B_0$: 
\begin{equation}\label{cov_rel3}
\BmZ{1}{0}\cover{\hF_*} B_0,  \quad T(\BmZ{1}{0}) \cover{\hF^{-1}_*} B_0.
\end{equation} 

\medskip

\noindent {\it Step 4.} Let $\tilde{\Bm{1}{0}} \Supset {\Bm{1}{0}} $ be any set such that $\eps \le {\rm dist}_{p \in \partial \tilde{\Bm{1}{0}} }(p,\Bm{1}{0}) \le 2 \eps$ for some sufficiently small $\eps>0$, chosen so that  $\hat{F}_*$ is analytic on $\tilde{\Bm{1}{0}}$. The set  $\BmZ{1}{1}=\hat{F}_*^{-1}\left(T(\tilde{\Bm{1}{0}})\right) \bigcap B_0$ satisfies 
\begin{equation}\label{cov_rel4}
\BmZ{1}{1} \cover{\hF_*} B_1, \quad T(\BmZ{1}{1}) \cover{\hF^{-1}_*} B_1.
\end{equation}


Set
$$\cC_{\hF_*}^\infty=\lim_{k \rightarrow \infty} \cC_{\hF_*}^k,$$
where 
$$\cC_{\hF_*}^k=\hF_*(\cC_{\hF_*}^{k-1}) \cap \hF_*^{-1}(\cC_{\hF_*}^{k-1}), \quad \cC_{\hF_*}^1=\hF_*(\cup_{i,j=0}^1 \!\!\!\phantom{a}_{1}B_i^j) \cap
\hF^{-1}_*(\cup_{i,j=0}^1 T(\!\!\!\phantom{a}_{1}B_i^j)).$$

We have proved the following:

\begin{lemma}
The map $\hF_*=F_*^{42}$ admits a locally maximal invariant set $\cC_{\hF_*}^\infty$ on which its action is homeomorphic to the full Bernoulli shift on $\{0,1\}^{\fZ}$.
\end{lemma}

\noindent {\it Proof of Main Theorem \ref{MTHM1} (2-4):}

The covering relations $(\ref{cov_rel1})$ imply that, for any $m$, there exist vertical sets $\Bm{m}{0}\subset \Bm{1}{0}$ and   $\Bm{m}{1}\subset \Bm{1}{1}$ such that $\hF_*^{m}$ is analytic on them, and $\Bm{m}{0} \cover{\hF_*^{m}} B_0$ and $\Bm{m}{1} \cover{\hF_*^{m}} B_1$.

Consider the map $\hF_*^{2^k}$ for some integer $k >1$.  The fixed point equation 
\begin{equation}\label{fpeq} 
\Lambda_*^{-k} \circ \hF^{2^k}_* \circ \Lambda_*^k=\hF_*
\end{equation}
implies that
\begin{eqnarray}
\nonumber \Lambda_*^{-k}(\Bm{n 2^k}{0}) &\cover{\hF_*^n} &\Lambda_*^{-k}(B_0) \equiv \Lambda_*^{-k-1}(B_1), \\
\nonumber  \Lambda_*^{-k}(T(\Bm{n 2^k}{0})) &\cover{\hF_*^{-n}} & \Lambda_*^{-k}(T(B_0)) \equiv \Lambda_*^{-k-1}(T(B_1)),\\
\nonumber \Lambda_*^{-k}(\Bm{n 2^k}{1}) &\cover{\hF_*^n} &\Lambda_*^{-k}(B_1),\\
\nonumber \Lambda_*^{-k}(T(\Bm{n 2^k}{1})) &\cover{\hF_*^{-n}}& \Lambda_*^{-k}(B_1).
\end{eqnarray}

Below we will use the notation that  $A \cover{\hF_*^n}  B$ if there exist some vertical set $\tilde{A} \subset A$ (depending on $n$ and $A$) on which $\hF_*^n$ is analytic, and $\tilde{A} \cover{\hF_*^n}  B$. $A\bcover{\hF_*^n}B$ and $A\gcover{\hF_*^n}B$ are used similarly. The covering relations $(\ref{cov_rel1})$ together with the fixed point equation $(\ref{fpeq})$ imply that, in this notation, we get for any integer $k >0$
\begin{eqnarray}
\nonumber\Lambda_*^{k}(B_1) \gcover{\hF_*^{2^k}}  \Lambda_*^{k-1}(B_1) \gcover{\hF_*^{2^{k-1}}}  \ldots \gcover{\hF_*^2} B_1,\\
\nonumber\Lambda_*^{k}(T(B_1)) \gcover{\hF_*^{-2^k}}  \Lambda_*^{k-1}(T(B_1)) \gcover{\hF_*^{-2^{k-1}}}  \ldots \gcover{\hF_*^{-2}} T(B_1).
\end{eqnarray}

We have, therefore, shown that for any two natural  $N$ and $M$, there exists an integer $n \equiv 2^N$, such that
\begin{eqnarray}
\nonumber
\Lambda^{N}_*(B_1) & \gcover{\hF_*^n}  \Lambda^{N-1}_*(B_1) \gcover{\hF_*^n}  \ldots \gcover{\hF_*^n} B_1 \gcover{\hF_*^n}   \Lambda^{-1}_*(B_1) \gcover{\hF_*^n} \\ 
\nonumber  \ldots & \gcover{\hF_*^n} \Lambda^{-M+1}_*(B_1) \gcover{\hF_*^n}  \Lambda^{-M}_*(B_1) 
 \gcover{\hF_*^{n}} \Lambda^{-M+1}_*(B_1) \gcover{\hF_*^{n}}\\ 
\nonumber \ldots & \gcover{\hF_*^{n}} \Lambda^{-1}_*(B_1) \gcover{\hF_*^{n}}  B_1 \gcover{\hF_*^{n}}  \ldots \gcover{\hF_*^{n}}
\Lambda^{N-1}_*(B_1) \gcover{\hF_*^n}  \Lambda^{N}_*(B_1)  
\end{eqnarray}

Set 
$$
\Delta_{N,M} = \bigcup_{k=-M}^{N} \Lambda^k_* \left(\Bm{\phi(k)}{0}\right)\cup\Lambda^k_* \left(\Bm{\phi(k)}{1}\right),
$$
where $\phi(k)=1$ if $k\geq 0$ and $\phi(k)=2^{N-k}$ if $k<0$, and define recursively:

$$ \cZ_{\hF^n_*}^1 \equiv \hF^n_*(\Delta_{N,M}) \cap \hF^{-n}_*(T(\Delta_{N,M})) \quad {\rm and} \quad \cZ_{\hF^n_*}^k \equiv  \hF^n_*(\cZ_{\hF^n_*}^{k-1}) \cap \hF^{-n}_* (\cZ_{\hF^n_*}^{k-1}).$$  

Clearly, the set $\cZ_{F^n_*}\equiv \cZ_{F^n_*}^\infty$ is a locally maximal invariant set for $F_*^n$, appearing in part $3)$ of the Main Theorem 1.



$F^n_* \arrowvert_{\cZ_{F^n_*}}$ is homeomorphic to a topological Markov chain, which is mixing. This together with ${\rm dist}(\Lambda_*^{-k}(B_1),0) \converge{{k \rightarrow \infty}} \infty$ and  ${\rm dist}(\Lambda_*^{k}(B_1),0) \converge{{k \rightarrow \infty}} 0$  implies the existence of {\it unbounded} and {\it oscillating} orbits. This proves part $4)$ the Main Theorem 1.


Finally, we notice that if $\cO_0$ is a heteroclinic orbit between points $z_1$ and $z_{0}$ for the iterate $\bar{F}_* \equiv F_*^4$ (or vice versa), then the fixed point equation $(\ref{fpeq})$ implies that, on one hand $\cO_k \equiv \Lambda_*^{k}(\cO_0)$ is a heteroclinic orbit for iterate $\bar{F}_*^{2^{k}}$ between periodic points $z_{k+1}$ and $z_k$:
$$z_k\equiv \Lambda^k_*(z_0);$$
on the other hand  if $p \in \cO_{-k} \equiv \Lambda_*^{-k}(\cO_0)$ then $\Lambda_*^k(p) \in \cO_0$, and so is $q=\bar{F}^{2^k}_*(\Lambda_*^k(p))$, and 
$$\bar{F}_*(p)=\Lambda_*^{-k}(\bar{F}^{2^k}_*(\Lambda_*^k(p)))=\Lambda_*^{-k}(q)  \in \cO_{-k}.$$

\begin{figure}[ph]
\vspace{2.5cm} a) \vspace{-4.5cm}
\begin{center}
\includegraphics[width=0.5\textwidth]{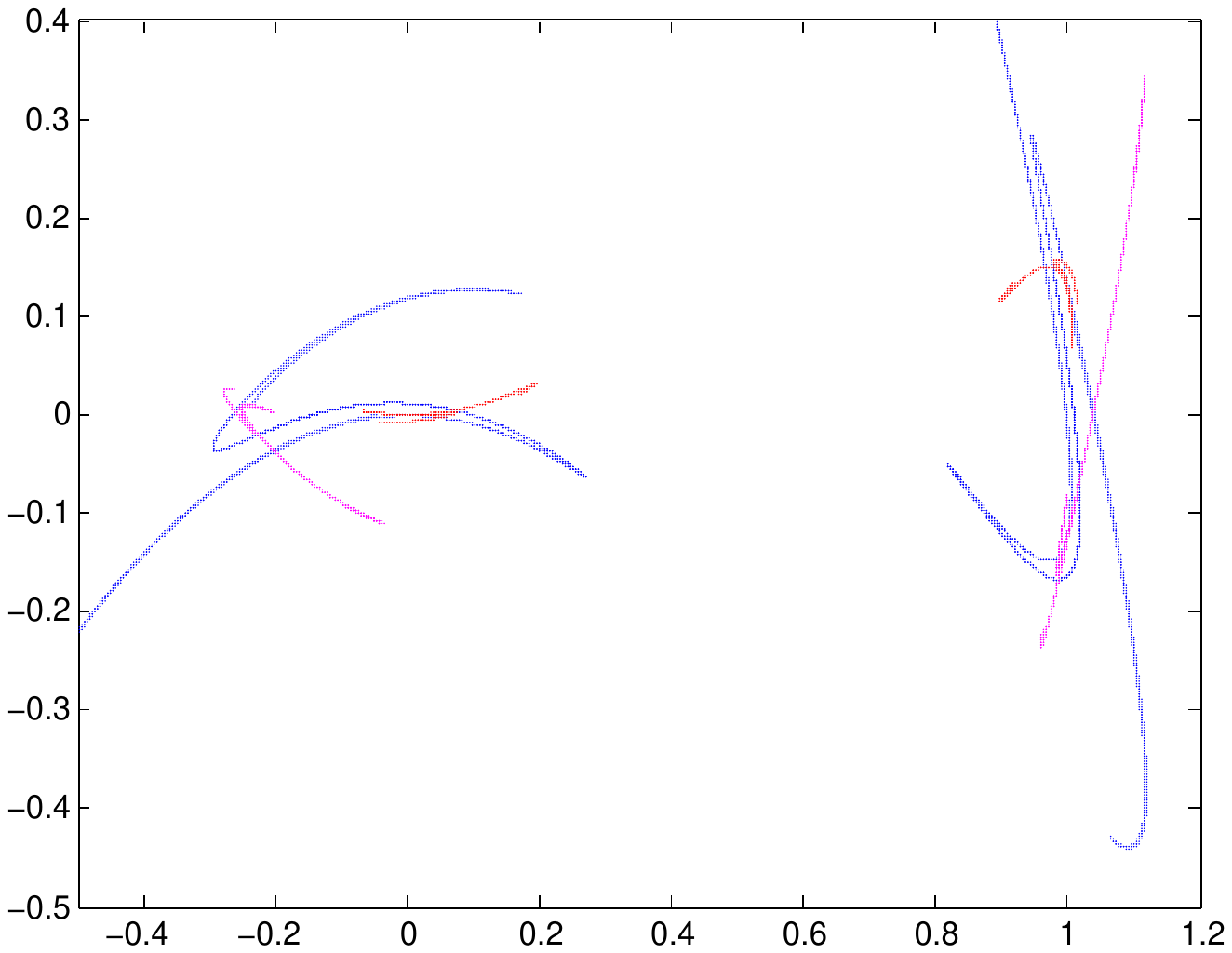} 
\vspace{-2cm}$$\uparrow \Lambda_*$$
\end{center}
\vspace{2cm} b) \vspace{-5.0cm}
\begin{center}
\includegraphics[width=0.5\textwidth]{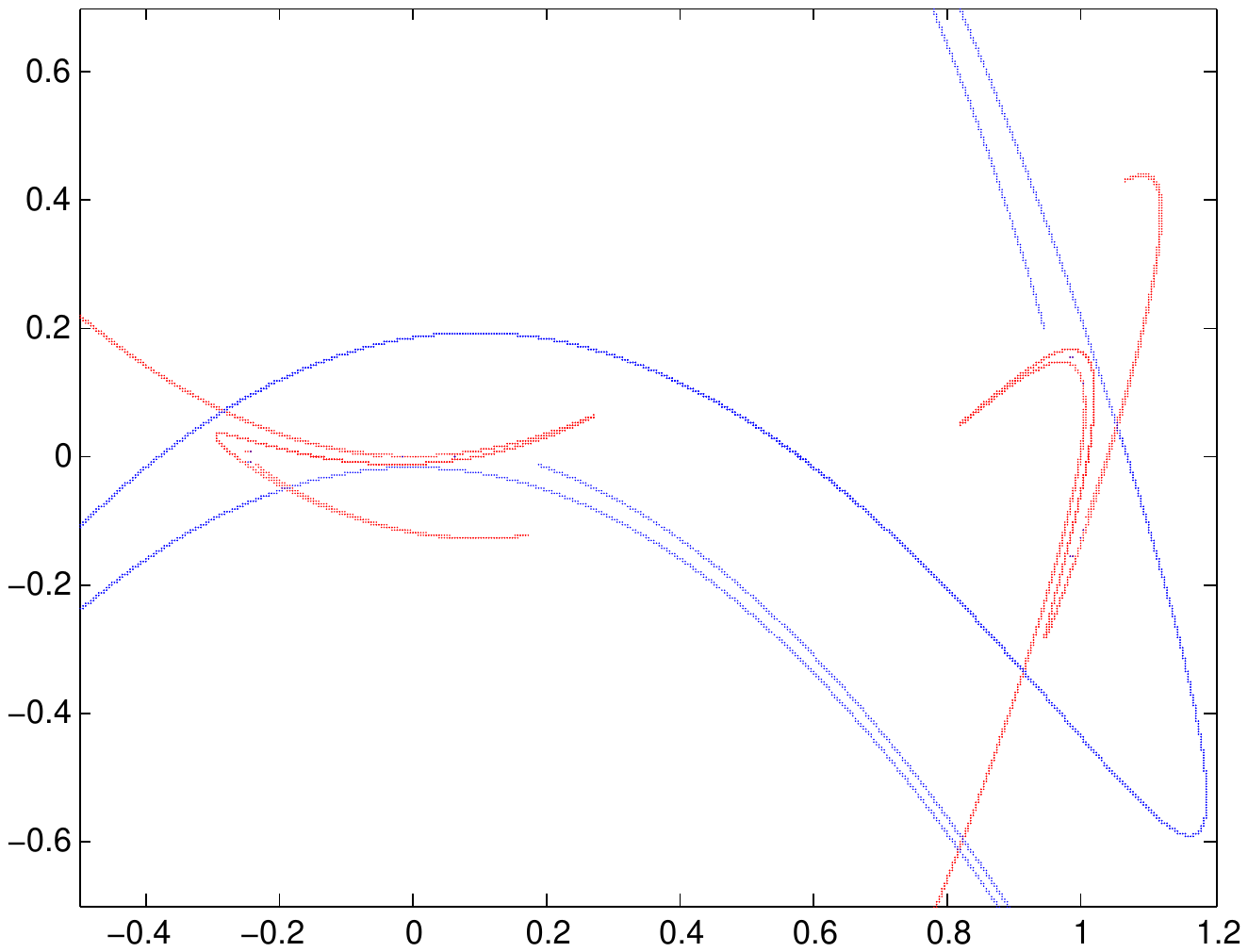} 
\vspace{-2cm}$$\downarrow \Lambda_*^{-1}$$
\end{center} 
\vspace{2.5cm} c) \vspace{-5.0cm}
\begin{center}
\includegraphics[width=0.5\textwidth]{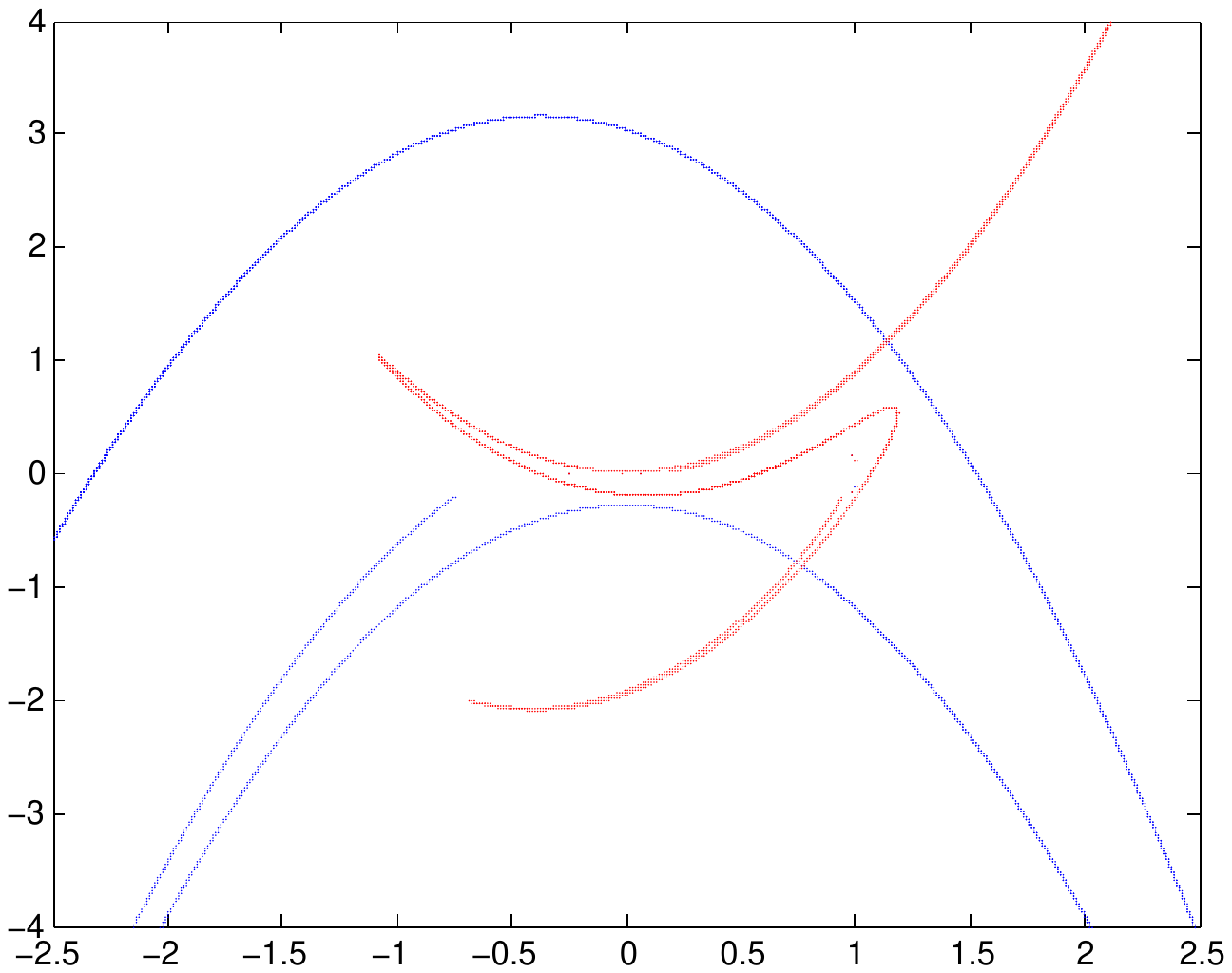} 
\caption{\it Intersections of a) the stable manifolds of $z_1$ and $F_*(z_1)$ (blue) with the unstable ones for $z_2$, $F_*(z_2)$, $F_*^2(z_2)$ and $F_*^3(z_2)$ (red and magenta); b) the stable manifold of $z_0$ (blue) with the unstable ones for $z_1$ and $F_*(z_1)$ (red); c) the stable manifold of $z_{-1}$ (blue)  with the unstable one for $z_0$ (red).}\label{TangleFig}
\end{center}
\end{figure}

Therefore, $\cO_{-k}$ contains a heteroclinic orbit between fixed points $z_{1-k}$ and $z_{-k}$ for iterate $\bar{F}_*$. Stable and unstable manifolds of points $z_k$ and $z_{k+1}$ intersect in a two-sided tangle, see Figure \ref{TangleFig}. This proves part $2)$ of the Main Theorem 1.

{\flushright {$\Box$}}

\section{Distortion tools}

We will say that a pair $(\cK,\psi)$ is a {\it dynamically  defined Cantor set} if $\cK \subseteq R$ is a Cantor set and $\psi: \cK \mapsto \cK$ is a locally Lipschitz expanding map, topologically conjugate to the full Bernoulli shift $\sigma: \{0,1\}^\fN \mapsto \{0,1\}^\fN$.

A {\it Markov partition}  $\cP =\{\cK_0,\cK_1\}$ of $(\cK,\psi)$ is a partition of $\cK$ in two disjoint subsets, $\cK=\cK_0 \cup \cK_1$, $\cK_0 \cap \cK_1 = \emptyset$, such that $\psi\arrowvert_{\cK_i}:\cK_i \mapsto \cK$ is a strictly monotone Lipschitz expanding homeomorphism.

Given a sequence $(a_0, \ldots, a_{n-1}  ) \in \{0,1\}^n$, denote
$$\cK(a_0,\ldots,a_{n-1}) =\bigcap_{i=0}^{n-1} \psi^{-i}(\cK_{a_i}).$$

A bounded component of $\fR \setminus \cK$ is called a ${\it gap}$ of $\cK$. Gaps can be ordered as follows. The unique gap of order $0$ is the interval $U_0 \equiv \hat{\cK} \setminus \left(\hat{\cK_0} \cup \hat{\cK_1}\right)$, where $\hat{A}$ denotes the convex hull of a set $A \subseteq \fR$. A connected component of 
$$\hat{\cK} \setminus \bigcup_{(a_0,\ldots,a_1) \in \{0,1\}^n} \widehat{\cK(a_0,\ldots,a_{n-1})}$$
is called a gap of order $n$ if it is not a gap of order $k < n$. Any gap of $\cK$ is a gap of some finite order.

\begin{definition}\label{dist_def}
Given a gap $U$ of $\cK$ of order $n$, we denote $L_U$ (respectively, $R_U$) the unique left (respectively, right) adjacent to $U$ interval of the form $\widehat{\cK(a_0,\ldots,a_{n-1})}$, $(a_0,\ldots,a_{n-1}) \in \{0,1\}^n$. 

The numbers
\begin{eqnarray}
\tau_L(\cK) &\equiv& \inf\left\{{|L_U| \over |U|}: {\rm U \ is \ a \ gap \  of \  \cK} \right\},\\
\tau_R(\cK) &\equiv& \inf\left\{{|R_U| \over |U|}: {\rm U \ is \ a \ gap \  of \  \cK} \right\} 
\end{eqnarray}
will be called the left and right thicknesses of $\cK$. The numbers
$$\tau_L(\cP) \equiv {|L_{U_0}| \over |U_0|  } \quad  {\rm and } \quad  
\tau_R(\cP) \equiv {|R_{U_0}| \over |U_0|  }$$ 
are called the left and right thicknesses of the Markov partition $\cP$.
Given a Lipschitz expanding map $\psi: I \mapsto \fR$, $I \subset \fR$, define the distortion of $\psi$ on $I$:
\begin{equation}
\hspace{-1.5cm}{\rm Dist}(\psi,I) \equiv {\rm sup}_{\{ x,y,z \in I: z\ne x, y\ne x \}} {\rm log} \left\{ { |\psi(y)-\psi(x)||z-x| \over |\psi(z)-\psi(x)||y-x| }  \right\}  \in [0,\infty]
\end{equation}
( notice, $\psi(z)\ne \psi(x)$, $\psi(y)\ne \psi(x)$). 
\end{definition}


\begin{definition}
The distortion of the dynamically defined Cantor set $(\cK,\psi)$ is
\begin{equation}
{\rm Dist}_{\psi}(\cK) \equiv {\rm sup}_{(a_0,\ldots,a_{n-1}) \in \{0,1\}^n} {\rm Dist}(\psi^n,\cK(a_0,\ldots,a_{n-1})).
\end{equation}
\end{definition}

\begin{lemma}\label{dist_lemma}
(see \cite{PT}, \cite{Duarte1}) Let $(\cK,\psi)$ be a dynamically defined Cantor set with a Markov partition $\cP$ and distortion ${\rm Dist}_\psi(\cK)=c$. Then
\begin{equation}
e^{-c} \tau_L(\cP) \le \tau_L(\cK) \le e^{c} \tau_L(\cP), \quad e^{-c} \tau_R(\cP) \le \tau_R(\cK) \le e^{c} \tau_R(\cP).
\end{equation}
\end{lemma}
This lemma is important since it allows one to estimate the thicknesses of a Cantor set in terms of thicknesses of its Markov partition.

\begin{definition}\label{Markovcomp} (see \cite{Duarte1}, \cite{GoKa}) Define $\cF$ to be the set of all maps $F: \Delta_0 \cup \Delta_1 \mapsto \fR^2$ such that:
\begin{itemize}
\item[1)] $\Delta_0$ and $\Delta_1$ are compact sets, diffeomorphic to rectangles, with non-empty interior;
\item[2)] $F$ is $C^2$ on a neighbourhood of $\Delta_0 \cup \Delta_1$, that maps  $\Delta_0 \cup \Delta_1$ diffeomorphically onto its image;
\item[3)] the locally maximal invariant set $\cC_F=\cap_{n \in \fZ} F^{-n}(\Delta_0 \cup \Delta_1)$ is a hyperbolic set, and the action of $F$ on $\cC_F$ is conjugated to the Bernoulli shift $\sigma: \{0,1\}^\fZ \mapsto \{0,1\}^\fZ$.
\item[4)] $\cP=\{\Delta_0,\Delta_1\}$ is a Markov partition for $F: \cC_F \mapsto \cC_F$, in particular, $F$ has two fixed points $p_0 \in \Delta_0$ and $p_1 \in \Delta_1$, whose stable and unstable manifolds contain the boundaries of $\Delta_0$ and $\Delta_1$. 
\end{itemize}
\end{definition}

We will now introduce our main tool for computing the Hausdorff dimension --- the Duarte Distortion Theorem (see \cite{Duarte1}, \cite{GoKa}).

\begin{definition}\label{classF}
Given  positive constants $C$, $\epsilon$ and $\gamma$ define $\cF(C,\epsilon,\gamma)$ to be the class of maps $F: \Delta_0 \cup \Delta_1 \mapsto \fR^2$, $F \in \cF$, such that 
\begin{itemize}
\item[1)] ${\rm diam}(\Delta_0 \cup \Delta_1) \le 1, \quad {\rm diam}(F(\Delta_0) \cup F(\Delta_1)) \le 1$;  
\item[2)]the derivative of $F$, $DF(x,u)=\left[a(x,u) \quad b(x,u) \atop c(x,u) \quad  d(x,u)\right]$, where $a$, $b$, $c$ and $d$ are $C^1$, satisfies on $\Delta_0 \cup \Delta_1$
\begin{itemize}
\item[a)]${\rm det} DF(x,u) \equiv 1$,
\item[b)]$|d| <1<|a| \le {C \over \epsilon}$, 
\item[c)]$|b|,|c| \le \epsilon (|a|-1)$;
\end{itemize}
\item[3)]functions $\tilde{a} \equiv a \circ F^{-1}$, $\tilde{b} \equiv b \circ F^{-1}$, $\tilde{c} \equiv c \circ F^{-1}$, $\tilde{c} \equiv c \circ F^{-1}$ on $F(\Delta_0) \cup F(\Delta_1)$, appearing in $DF^{-1}(x,u) \equiv \left[\phantom{-}\tilde{d} \quad -\tilde{b} \atop -\tilde{c} \quad \phantom{-}\tilde{a} \right] $, satisfy
\begin{eqnarray}
\label{gamma1} &a)& \quad \left| {\partial \tilde{b} \over  \partial x }\right|= \left| {\partial \tilde{d} \over  \partial u }\right|, \quad   \left| {\partial \tilde{b} \over  \partial u }\right|, \quad \left| { \partial \tilde{c} \over  \partial x } \right|, \quad  \left| { \partial \tilde{a} \over  \partial x }\right| = \left| {\partial \tilde{c} \over  \partial y }\right| \le \gamma (|\tilde{a} |-1),\\
\label{gamma2}  &b)& \quad \left| {\partial {a} \over  \partial u }\right|= \left| {\partial {b} \over  \partial x }\right|, \quad   \left| {\partial {b} \over  \partial u }\right|, \quad \left| { \partial {c} \over  \partial x } \right|, \quad  \left| { \partial \tilde{c} \over  \partial u }\right| = \left| {\partial {d} \over  \partial x }\right| \le \gamma (|\tilde{a} |-1),\\
\label{gamma3}  &c)& \quad \left| {\partial \tilde{a} \over  \partial u }\right|, \quad \left| {\partial \tilde{d} \over  \partial x }\right| \le \gamma |\tilde{a}| (|\tilde{a}|-1),\\
\label{gamma4} &d)& \quad  \left| {\partial {a} \over  \partial x }\right|, \quad \left| {\partial {d} \over  \partial u }\right| \le \gamma |\tilde{a}| (|\tilde{a}|-1);
\end{eqnarray} 
\item[4)] the variation of $\log|a(x,u)|$ in each rectangle $\Delta_i$ is less or equal to $\gamma(1-\alpha_i)$, $\alpha_i=\max_{(x,u)\in \Delta_i} |a(x,u)|$;
\item[5)] the gap sizes satisfy:
$${\rm dist}(\Delta_0,\Delta_1) \ge { \epsilon \over \gamma}, \quad {\rm dist}(F(\Delta_0),F(\Delta_1)) \ge { \epsilon \over \gamma}.$$
\end{itemize}
\end{definition}

\begin{remark}
Unlike in  \cite{Duarte1} and \cite{GoKa} we do not require the eigenvalues of the maps $F \in \cF$ at the fixed points $p_0$ and $p_1$ to be positive. The positivity of the eigenvalues is important only in the construction of the boundaries of the two components $\Delta_0$ and $\Delta_1$ of the Markov partition $\cP$ (see Lemma $\ref{lDeltaTilde}$ for the construction of $\cP$ in our case), but does not enter the proof of  Theorem \ref{Duarte}.   
\end{remark}

We will now continue with the notion of stable and unstable Cantor sets. Denote the stable and unstable foliations of $\cC_F$ as $\cF^s$ and $\cF^u$:
\begin{eqnarray}
\label{Fs} \cF^s &\equiv &\{{\rm connected \,\, comp. \,\, of} \,\, \cW^s(\cC_F) \cap \left( \Delta_0' \cup \Delta_1'\right)\},\\
\label{Fu} \cF^u &\equiv& \{{\rm connected \,\, comp. \,\, of} \,\, \cW^u(\cC_F) \cap \left( F(\Delta_0') \cup F(\Delta_1') \right)\}.
\end{eqnarray}
Also, define
$$\cI^s \equiv \cW^s_{\loc}(p_0) \cap \Delta_0' \quad {\rm and} \quad  \cI^u \equiv \cW^u_{\loc}(p_0) \cap \Delta_0.$$
Then, the Cantor sets 
\begin{equation}
\cK^s \equiv \cC_F \cap \cI^u  \quad {\rm and} \quad
\cK^u \equiv \cC_F \cap \cI^s 
\end{equation}
can be identified with the set of stable leaves of $\cF^s$, respectively unstable leaves of $\cF^u$. Define the projections $\pi_s: \cC_F \mapsto \cK^s$ and $\pi_u: \cC_F \mapsto \cK^u$ in the following way: if $p \in \cC_F$ then
$$\pi_s(p)=\cW^s_{\loc}(p) \cap \cI^u \quad {\rm and} \quad \pi_u(p)=\cW^u_{\loc}(p) \cap \cI^s,$$
and the maps $\psi^s: \cK^s \mapsto \cK^s$ and $\psi^u: \cK^u \mapsto \cK^u$ as
$$\psi^s \equiv \pi_s \circ F \quad {\rm and} \quad \psi^u \equiv \pi_u \circ F^{-1}.$$

The pairs $(\cK^s,\psi^s)$ and  $(\cK^u,\psi^u)$ are dynamically defined Cantor sets with Markov partitions 
$$\cP^s \equiv \{ \cI^u \cap \Delta_0, \cI^u \cap \Delta_1 \} \quad {\rm and} \quad  \cP^u \equiv \{ \cI^s \cap F(\Delta_0), \cI^s \cap F(\Delta_1) \}.$$

\begin{theorem} \label{Duarte}({\it Duarte Distortion Theorem, \cite{Duarte1}}) If $F \in \cF(C,\epsilon,\gamma)$, then the distortion of the dynamically defined Cantor sets $(\cK^{u,s},\psi_{u,s})$, $\psi_{u,s}=\pi_{u,s} \circ F$, is bounded:
\begin{equation}
{\rm Dist}_{\psi_{u,s}}(\cK^{u,s}) \le  D(C,\epsilon,\gamma) \equiv 4(C+3)\gamma +2 \epsilon.
\end{equation} 
In particular,
\begin{eqnarray}
\nonumber e^{-D(C,\epsilon,\gamma)} &\tau_L(\cP^{u,s})& \le \tau_L(\cK^{u,s}) \le  e^{D(C,\epsilon,\gamma)} \tau_L(\cP^{u,s}),\\
\nonumber e^{-D(C,\epsilon,\gamma)} &\tau_R(\cP^{u,s})& \le \tau_R(\cK^{u,s}) \le  e^{D(C,\epsilon,\gamma)} \tau_R(\cP^{u,s}).
\end{eqnarray}
\end{theorem}

The Duarte Distortion Theorem has been used as a powerful tool in conservative dynamics  on several instances. It has been proved, and applied in \cite{Duarte1} to demonstrate that accumulation of a locally maximal invariant set by periodic elliptic points is generic in two parameter families. Furthermore, \cite{GoKa} uses it to show that generic unfoldings of homoclinic tangencies of two-dimensional area-preserving diffeomorphisms give rise to hyperbolic sets of dimension arbitrarily close to $2$. To show this, the authors of \cite{GoKa} complement the Duarte Distortion theorem with the following results, which will be useful to us:

\begin{lemma}\label{thicknesslemma}
The Hausdorff dimension of a Cantor set $\cK$ satisfies 
$${\rm dim}_H \cK > d,$$ 
where $d$ is the solution of 
$$\tau_L(\cK)^d+\tau_R(\cK)^d=(1+\tau_L(\cK)+\tau_R(\cK))^d.$$
\end{lemma}    

\medskip

\begin{lemma}
Suppose that for given $t_L>0$, $t_R>0$ the solution $d$ of the equation
$$t_L^d+t_R^d=(1+t_L+t_R)^d$$
is in $(0,1)$, then
\begin{equation}
d > \max \left\{ {\log\left(1+{t_R \over 1+t_L} \right) \over \log\left(1+{1+t_R \over t_L} \right) },   {\log\left(1+{t_L \over 1+t_R} \right) \over \log\left(1+{1+t_L \over t_R} \right) }    \right\}.
\end{equation}
\end{lemma}

\section{A horseshoe for $\bG$}\label{topHorseShoe}
In this Section we demonstrate the existence of a horseshoe for the third iterate of $G\in\bG$. We start with definitions.

\begin{definition}\label{full_comp} (\underline{Full component})
Suppose $\Delta \subset \cU \subset \fR^n=\fR^k \oplus \fR^l$ is homeomorphic to a rectangle (a set of the form $D_1 \times D_2 \subset  \fR^k \oplus \fR^l$), and let $F: \cU \mapsto \fR^n$ be a diffeomorphism. A connected component $\Delta_0=F(\bar{\Delta}_0)$ of $\Delta \cap F(\Delta)$ is called {\it full}, if 
\begin{itemize}
\item[1)]$\tilde{\cP}_2 (\bar{\Delta}_0)=D_2$,
\item[2)]for any $z \in \bar{\Delta}_0$, $\tilde{\cP}_1\arrowvert_{F(\bar{\Delta}_0 \cap (D_1 \times \tilde{\cP}_2(z) )  )  }
$ is a bijection onto $D_1$.
\end{itemize}
Here, $\tilde{\cP}_{1,2} \equiv \cP_{1,2} \circ h$, $h$ the homeomorphism $\Delta_0 \, {\approx \atop \mbox{\it h}} \, D_1 \times D_2$, $\cP_{1}$ and $\cP_2$  are the canonical projections on $\fR^k$ and $\fR^l$. 
\end{definition}

\medskip

\begin{definition}\label{horseshoe_def}(\underline{Two-component horseshoe})
Let $\cU \subset \fR^n$ be an open set, then the set $\Delta \subset \cU$ homeomorphic to a rectangle $D_1 \times D_2$  is called a {\it two-component horseshoe} for the diffeomorphism $F: \cU \mapsto \fR^n$ if $\Delta \cap F(\Delta)$ contains at least two full components $\Delta_0$ and $\Delta_1$ such that,
\begin{itemize}
\item[1)] $\tilde{\cP}_2(\Delta_0 \cup \Delta_1) \subset {\rm int}\,D_2$,  $\tilde{\cP}_1(F^{-1}(\Delta_0 \cup \Delta_1)) \subset {\rm int}\, D_1$;
\item[2)] $D F\arrowvert_{F^{-1}(\Delta_0 \cup \Delta_1)} $ preserves and expands an unstable (``horizontal'') cone family on $F^{-1}(\Delta_0 \cup \Delta_1)$;
\item[2)] $D F^{-1} \arrowvert_{\Delta_0 \cup \Delta_1} $ preserves and expands a stable (``vertical'') cone family on $\Delta_0 \cup \Delta_1$.
\end{itemize}
\end{definition} 

\medskip
By a standard result of the theory of dynamical systems, if $\Delta$ is a two-component horseshoe for $F$, then the action of $F$ on the locally maximal hyperbolic set $C_F \equiv \cap_{i=-\infty}^{\infty} F^i(\Delta)$ is homeomorphic to the Bernoulli shift $\sigma_2$ on $\{0,1\}^{\fZ}$.

To demonstrate the existence of a horseshoe for the third iterate of  $G \in \bG$, we construct two rectangles: $\Delta'_0$ which contains 
the fixed point $p_0$ and $\Delta_1'$ which contains the period three point $p_1$. We will prove that these two rectangles constitute two full components of a horseshoe, and are such that $G(\Delta_0') \cap \cW^s_\loc(p_0)$ and $G(\Delta_1') \cap \cW^s_\loc(p_0)$ contain the components of the {\it Markov partition} of the locally maximal hyperbolic set $\cK^u_{*}$.

\begin{table}[b]
\begin{center}
\begin{footnotesize}
\begin{tabular}{c|ccc}
Component & Centre & ``Stable'' Scale & ``Unstable'' Scale\\ \hline
$\Delta_0'$ & $(0.670198,  0.0)$ & $0.083$ & $0.083$ \\
$\Delta_1'$ & $(-0.441811, 0.0)$ & $0.0655$ & $0.0655$
\end{tabular}
\begin{center}
\caption{The rectangles that approximate the Markov partition for the horseshoe of $\bG$. The spanning vectors are approximate stable and unstable vectors for $D\bG(\bp_{0,1})$. The ``stable'' vectors for $\bp_0$ and $\bp_1$, respectively, are ${\bf e}^s_0=(0.788578889012330, -0.614933602760558)$ and ${\bf e}^s_1=(0.750925931392967773, 0.660386436536671957)$. The ``unstable'' vectors are ${\bf e}^u_0=T({\bf e}^s_0)$ and  ${\bf e}^u_1=T({\bf e}^s_1)$.  The length of the sides of the rectangles $\Delta_0'$ and $\Delta_1'$ is $2 \cdot {\rm stable/unstable \quad \!\! scale} \cdot |{\bf e}^{u,s}_{0,1}|$.  The covering sequence, which holds for $\bG$, is 
$\Delta_0'\Rightarrow \Delta_0'\Rightarrow \Delta_1'\Rightarrow \Delta_1'\Rightarrow \Delta_0'$.}\label{bigHorseshoe}
\end{center}
\end{footnotesize}
\end{center}
\end{table}

The  centres of the rectangles and sizes of the spanning vectors,  are given in Table \ref{bigHorseshoe}. Note that these boxes are  symmetric with respect  to the involution $T: (x,u) \mapsto (x,-u)$. Therefore, we only have to prove that the covering relations hold for $\bG$.

\begin{theorem}\label{tP3TopHorse}
For any $G\in\bG$, the inverse map $G^{-1}$ admits a horseshoe whose two full components are the sets $\Delta_0'$ and $\Delta_1'$, defined in Table \ref{bigHorseshoe}. $\Delta_0'$ and $\Delta_1'$ are independent of $G$.
\end{theorem}

\noindent {\it Proof.} To demonstrate the existence of the horseshoe we first use the routine {\tt checkCoveringRelations} to show that
$$
{\Delta}_0' \cover{\bG} {\Delta}_0' \cover{\bG}{\Delta}_1' \cover{\bG}{\Delta}_1' \cover{\bG}{\Delta}_0',
$$
where the $h$-sets are as in Table \ref{bigHorseshoe}.
This implies that $1)$ in the Definition $\ref{horseshoe_def}$ holds.
 
The remainder of the proof of the existence of a horseshoe is done for the map $\tbG^{-1} \equiv \cT \circ \bG^{-1} \circ \cT^{-1}$, where $\cT(x,u)$ is some coordinate transformation that approximately diagonalizes the derivative $D\bG$ on $\Delta_0' \cup \Delta_1'$. (See \cite{GP} for programs used in this part of the proof).  We have used the following coordinate change: 
\begin{eqnarray}\label{linChange}
\!\!\!\!\!\!\!\!\!\!\!\!\!\!\!\!\!\!\!\!\! \nonumber \cT(x,u) = 0.55 \, (x-0.5, u&+&0.0265461632116977382 \\
\!\!\!\!\!\!\!\!\!\!\!\!\!\!\!\!\!\!\!\!\!\nonumber     &-&0.658388496175704694(x-0.577619) \\
\!\!\!\!\!\!\!\!\!\!\!\!\!\!\!\!\!\!\!\!\!\nonumber     &-&0.611583723237529069(x-0.577619)^2 \\
\!\!\!\!\!\!\!\!\!\!\!\!\!\!\!\!\!\!\!\!\!    &+&0.102408008658008658(x-0.577619)^3).
\end{eqnarray} 

Denote $\tilde{\Delta}_{0,1}'=\cT(\Delta_{0,1}')$.  We show that there exist an invariant unstable cone field $\tilde{\cC}_u$ on $\tilde{\Delta}_0' \cup \tilde{\Delta}_1'$,
\begin{eqnarray}
\label{tCu1} \tilde{\cC}_u(p) &\equiv& \{{\bf v}=(v_1,v_2) \in \fR^2: -0.15 \le {v_2 \over v_1} \le 0.2 \}, \quad p \in \tilde{\Delta}_0',\\
\label{tCu2} \tilde{\cC}_u(p) &\equiv& \{{\bf v}=(v_1,v_2) \in \fR^2: -0.1 \le {v_2 \over v_1} \le 0.05\}, \quad p \in \tilde{\Delta}_1'.
\end{eqnarray} 

We have verified that
\begin{equation}
\!\!\!\!\!\!\!\!\!\!\!\!\!\!\!\!\!\!\!\!\!\!\!\!    \!\!\!\!\!\!\!\!\!\!\!\!\!\!\!\!\!\!\!\!\!\label{consGu}  D\tbG(p) \cdot {\bf v} \in {\rm int} \, \tilde{\cC}_u(\tbG(p)), \, {\rm for \, all} \quad {\bf v} \in \tilde{\cC}_u(p), \quad p \in \tilde{\Delta}_0' \cup \tilde{\Delta}_1',
\end{equation}
and that vectors inside the cone field are expanded:
$$ | D\tbG(p) \cdot {\bf v} | =A | {\bf v}|, \quad {\bf v} \in {\rm int} \, \tilde{\cC}_u(\tbG(p)),$$
where the expansion rate $A$ is in $\bA=[A_-,A_+]$
\begin{equation}\label{Apm}
A_-= 6.091, \quad A_+=26.892.
\end{equation}

Verification of invariance and expansion of the cone fields is carried out (among other things) in the routine {\tt Distortion} of \cite{GP}. 

Reversibility and symmetry of the partition imply that there exists a stable cone field 
\begin{equation}
\label{tCs}  \tilde{\cC}_s(p) \equiv D \cT(\cT^{-1}(p)) (T(D \cT^{-1}(p) \cC_u(p)))
\end{equation} 
on $\tbG(\tilde{\Delta}_0' \cup \tilde{\Delta}_1')\cap \left( \tilde{\Delta}_0' \cup \tilde{\Delta}_1' \right)$, which is mapped into its interior and expanded with the rate $A$ by  $D \tilde{G}^{-1}$.
\newline
{\flushright {$\Box$}}


\begin{lemma}\label{transversal_cones}
The cone fields 
\begin{eqnarray}
\label{Cu} \cC_u(p)&=&D \cT^{-1}(\cT(p)) \, \tilde{\cC}_u(\cT(p)), \quad p \in \Delta_0' \cup \Delta_1',\\
\label{Cs} \cC_s(q)&=&T(\cC_u(q)), \quad q \in \bG(\Delta_0' \cup \Delta_1')\cap \left(\Delta_0' \cup \Delta_1'\right),
\end{eqnarray}
are transversal in the sense that the angle between any ${\bf u} \in \cC_u(q)$ and  ${\bf v} \in \cC_s(q)$, $q \in  \bG(\Delta_0' \cup \Delta_1')\cap \left(\Delta_0' \cup \Delta_1'\right)$ is bounded from below. 

Furthermore, the leaves of the foliations $\cF^s$ and $\cF^u$, defined in $(\ref{Fs})$ ---$(\ref{Fu})$, are graphs over the $x$-axis.
\end{lemma}
\noindent {\it Proof.} We verify on the computer (this is done in the routine {\tt Distortion} in \cite{GP}) that any ${\bf v} \in  
D \cT^{-1}(\cT(p)) \tilde{\cC}_u(\cT(p))$, $p \in \Delta_0'$, satisfies:
$${ v_2 \over  v_1 } \ge  0.461,$$
and  that any ${\bf v} \in  D \cT^{-1}(\cT(p)) \tilde{\cC}_u(\cT(p))$, $p \in \Delta_1'$, satisfies:
$${ v_2\over   v_1 }\le -0.675.$$
This implies that the angle between any ${\bf u} \in \cC_u(q)$ and  ${\bf v} \in \cC_s(q)$, $q \in  \bG(\Delta_0' \cup \Delta_1')\cap \left(\Delta_0' \cup \Delta_1'\right)$ is bounded from below. 
\newline
{\flushright {$\Box$}}



\begin{figure}[h]
\vspace{4.0cm} a) \vspace{-4.0cm}
\begin{center}
\includegraphics[width=0.7\textwidth]{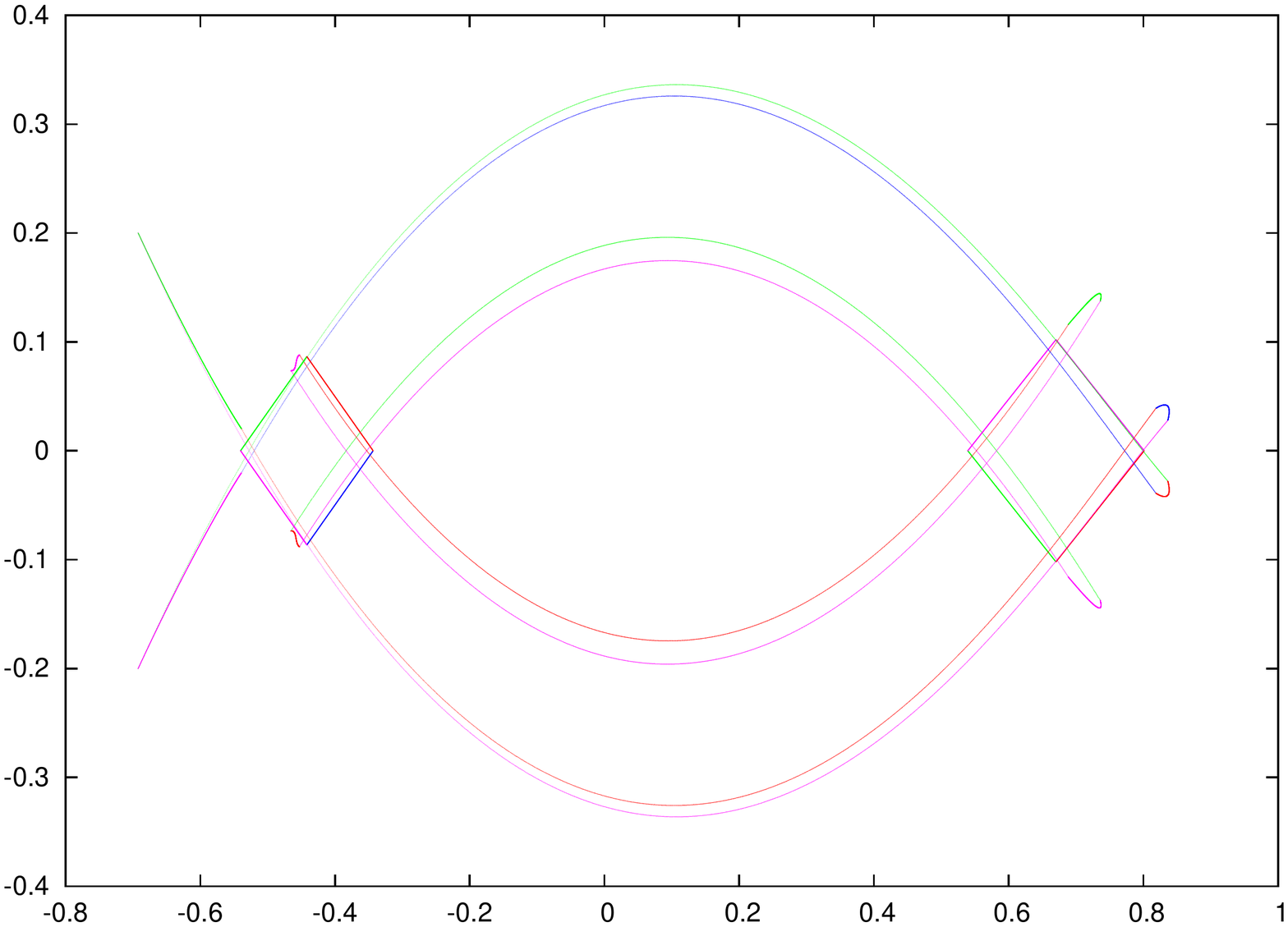} 
\end{center}
\vspace{3.5cm} b) \vspace{-4.0cm}
\begin{center}
\includegraphics[width=0.7\textwidth]{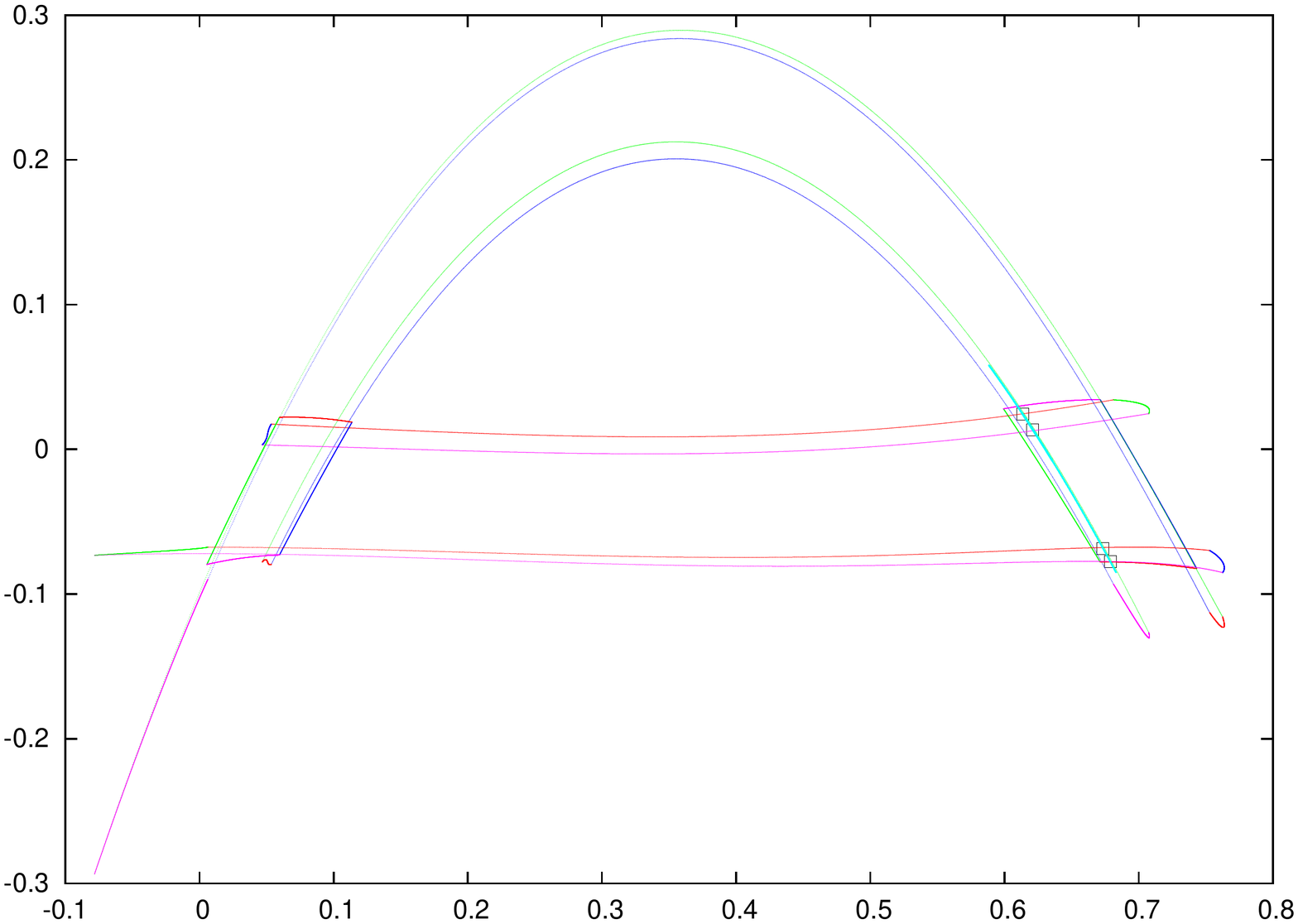} 
\end{center}
\begin{center}
\caption{a) The two components $\Delta_0'$ and $\Delta_1'$ of the horseshoe map across $\Delta_0' \cup \Delta_1'$ under $G$ and $G^{-1}$. $\Delta_0'$ is the diamond on the right, $\Delta_0'$ --- the diamond on the left. The entrance sets are given in red and magenta, the exit sets are in green and blue. The two sets whose longer sides in red and magenta are $G(\Delta_i')$, $i=0,1$, the sets with longer blue and green sides are $G^{-1}(\Delta_i')$, $i=0,1$. b) The same horseshoe in the coordinate $\cT$. The local stable manifold is given in cyan, the boxes signify the endpoints of the Markov partition for $\cK^u_*$. }\label{GHorsePic}
\end{center}
\end{figure}

The horseshoe for $G$ is illustrated in Figure \ref{GHorsePic}, for both original (a) and $\cT$ (b) coordinates.

Note that the full components $\Delta'_0$, and $\Delta'_1$ are chosen uniformly for all $G\in\bG$. This is not compatible with the definition of full components in Definition \ref{Markovcomp}. We can, however, prove that any $G\in\bG$ has full components, as in Definition \ref{Markovcomp}, contained in the components constructed in Theorem \ref{tP3TopHorse}. 

 \begin{lemma}\label{lDeltaTilde}
For any $G\in\bG$ there exist full components $\tD_0 \subset \tD^{'}_0$ and $\tD_1 \subset \tD^{'}_1$ as in Definition \ref{Markovcomp}.
\end{lemma}
\noindent {\it Proof.} (see \cite{GP} for programs) Let $G\in\bG$, by Theorem \ref{tP3TopHorse} 
$$\tilde p_1 \equiv \cT(p_1) \Subset  \tD^{'}_1\cap \tG(\tD^{'}),$$
hence there exists $\eps > 0$ such that $B_\eps(\tilde{p}_1) \subset \tD^{'}$. Therefore, $\tilde \cW^u_\eps \equiv \tilde \cW^u_\loc(\tilde{p}_1) \cap B_\eps(\tilde p_1) \subset \tD^{'}$. Let $t \mapsto \tilde \cW^u_\eps(t)$ be some parametrization of $\tilde \cW^u_\eps$. For sufficiently small $\eps$ the tangent vector $\tau(t)$  to  $\tilde \cW^u_\eps(t)$ is close to ${\bf e}^u_1(\tilde{p}_1)$, that is, there exists $\eps>0$ such that $\tau(t) \in \tilde{\cC}_u(\tilde \cW^u_\eps(t))$. 

The cone conditions imply that $\tilde \cW^u_\eps(t)$ is expanded in the forward (unstable) cone into ${\rm cvx}\,\tG(\tD^{'})$: there exists a positive integer $N$ such that $\tG^N(\tilde \cW^u_\eps) \cap (\tD_{0,1}^{'})^-_1 \ne \emptyset$ and  $\tG^N(\tilde \cW^u_\eps) \cap (\tD_{1,0}^{'})^-_2 \ne \emptyset$ where $(\tD_{0,1}^{'})^-_i$, $i=1,2$, denotes the two components of the exit set of $\tD_{0,1}^{'}$ (see Definition $\ref{hset}$). 

Set $b_0 \equiv \tG^N(\tilde \cW^u_\eps) \cap \tD_0^{'}$ and $b_1 \equiv \tG^N(\tilde \cW^u_\eps) \cap\tD_1^{'}$. $\tG(b_0)$, stretches across $\tD^{'}$, by the set of covering relations in Table \ref{bigHorseshoe}, that is  $\tG(b_0) \cap (\tD_{0,1}^{'})^-_1 \ne \emptyset$ and  $\tG(b_0) \cap (\tD_{1,0}^{'})^-_2 \ne \emptyset$.  We set $t_0 \equiv \tG(b_0) \cap \tD_0^{'}$ and $t_1 \equiv \tG(b_0)\cap\tD_1^{'}$. 

Via a similar construction for $\tG^{-1}$, we can show that there is a pair of stable leaves $l_1$ and $r_1$ such that $l_1 \cap (\tD_1^{'})^+_i \ne \emptyset$, $i=1,2$, and  $r_1 \cap (\tD_1^{'})^+_i \ne \emptyset$, $i=1,2$. Here $(\tD_{0,1}^{'})^+_i$, $i=1,2$, denote the two components of the entrance set of $\tD_{0,1}^{'}$ (see Definition $\ref{hset}$). Furthermore, there is a pair $l_0$ and $r_0$ with similar properties. 

Since $b_1$ and $t_1$ stretch across $\tD_1^{'}$ ``horizontally'', and $l_1$ and $r_1$ stretch across $\tD_1^{'}$ ``vertically'', $b_1$ intersects $l_1$ and $r_1$ at least once, and  $t_1$ intersects $l_1$ and $r_1$ at least once.   Clearly, $\tilde{p}_1$ is one of these intersection points. Therefore, $b_1$, $t_1$, $r_1$ and $l_1$ bound an open set $\cU \in \tD_1^{'}$ with at least one connected component. We define $\tilde{\Delta}_1$ as the closure of the connected component that contains $\tilde{p}_1$ in its boundary. 

$\tD_0$ is constructed in a similar way: it is the closure of the connected component bounded by $b_0$, $t_0$, $l_0$ and $r_0$ that contains $p_0$ in its interior.
\newline
{\flushright {$\Box$}}


\section{Hausdorff dimension of the locally maximal invariant set for  $\bG$}\label{HDC}

According to Lemma $\ref{thicknesslemma}$ and the Duarte Distortion Theorem $\ref{Duarte}$, the knowledge of the left and right thicknesses of the Markov partition for the unstable Cantor set $\cK^u$ together with the distortion of this Cantor set is sufficient to compute a lower bound on its Hausdorff dimension. 

\begin{lemma}\label{lKuEnc} 
For any $G\in\bG$ the intersections of the local stable manifold $\tilde \cW^s_\loc(\tilde p_0)$ with $\partial \tG(\tD_0)$ and $\partial \tG(\tD_1)$ are at the points $k_0^b, k_0^t$ and $k_1^b, k_1^t$, respectively.

Where,
\begin{eqnarray}
\nonumber k_0^b &\in& [0.041476215,0.04509142] \times [0.01162711,0.015242315], \\
\nonumber k_0^t &\in& [0.034578335,0.03698849] \times [0.023247785, 0.02565794],
\end{eqnarray}
and 
\begin{eqnarray}
\nonumber k_1^b &\in& [0.099814,0.1022241] \times [-0.07865165,-0.07624155], \\
\nonumber k_1^t &\in& [0.09487115,0.09660585] \times [-0.06980765,-0.06710275],
\end{eqnarray}
respectively.
\end{lemma}
\noindent {\it Proof.} (see \cite{JP} for programs)  We prove that there exists two points $q_1, q_2 \in \cW^u_\loc(p_1)\cap G(\Delta_1)\cap \cW^s(p_0)$, such that 
$$
G(q_1)\cap \cW^s_\loc(p_0) \neq \emptyset
$$  
$$
G(q_2) \cap \Delta_0 \neq \emptyset
$$
$$
G^2(q_2) \cap \Delta_1 \neq \emptyset
$$
$$
G^3(q_2) \cap \cW^s_\loc(p_0) \neq \emptyset 
$$
We now describe the construction of $q_1$. Numerical experiments indicate that there exists a heteroclinic orbit 
between  $p_0$, and the period three point $p_1$, approximately 
located at $(-0.527156,0)$. To prove this we construct a covering 
sequence of $h$-sets, as in Table \ref{tHeteroTraj}, consisting of: one 
box around $p_0$,  one box around $p_1$ and five intermediate 
boxes on the orbit; two of them close to the unstable vector of 
$p_1$ and three of them close to the stable vector of $p_0$.  We use the programs {\tt checkCoveringRelations} and {\tt checkConeConditions} to verify that the hypotheses of Theorem \ref{tCovRel} and Lemma \ref{lQF} are satisfied on 
the sequence of $h$-sets.

We emphasize that the choice of boxes is very delicate. To construct the 
sequence we first construct relatively large boxes around $p_0$ and 
$p_1$, spanned by their unstable and stable eigenvectors. 
Second, we reduce the size of these boxes until we are able to prove that the hypothesis of Lemma \ref{lQF} is satisfied. Third, we non-uniformly 
distribute points in the unstable eigendirection inside of the box containing the period 
three point; on the smallest scale the points are only separated by 
$10^{-11}$. Fourth, we iterate, non-rigorously, these points to locate a 
heteroclinic pseudoorbit, which starts within the box containing $p_1$ and ends in the box containing $p_0$. The 
first point on the orbit is in the unstable eigendirection of the period 
three point, at a distance of $0.000158885$ from our approximation of 
the period three point, and the last point is at a distance of 
$0.00177401$ from the fixed point. Finally, we use the points on this 
pseudoorbit as centres of our covering sequence. The first two intermediate 
points are on the local unstable manifold of the period three point, and 
we use its unstable and stable eigenvectors to span the corresponding 
$h$-sets. The other three intermediate points are on the local stable 
manifold of the fixed point, and we use its unstable and stable 
eigenvectors to span the corresponding $h$-sets. Note that if we change 
the location of the original point on the stable vector by as little as 
$10^{-8}$, the proof fails. The existence of $\tilde q_1$ follows from the set of covering relations constructed in Table \ref{tHeteroTraj}. 

To find $q_2$ we do a similar non-rigorous computation to locate a pseudoorbit with the properties above. The corresponding sequence of $h$-sets and covering relations is given in Table \ref{HeteroTraj2}. Finally, we put $k_1^b=\cT G(q_1)$, $k_0^t=\cT G^2(q_1)$, $k^b_0=\cT G^3(q_1)$, and $k_1^t=\cT G^3(q_2)$. By the construction of the components of the Markov partition $\tilde{\Delta}_0$ and $\tilde{\Delta}_1$, it is immediate that $\{k_0^t,k_0^b\}= \tilde \cW^s_\loc(\tilde p_0) \cap \partial \tG(\tD_0)$ and  $\{k_1^t,k_1^b\}= \tilde \cW^s_\loc(\tilde p_0) \cap \partial \tG(\tD_1)$. Evaluating $k^{b/t}_{0/1}$ at $\bq_{0/1}$ with $\bG$ yields the enclosures given in the statement of the lemma.

\begin{table}[b]
\begin{center}
\begin{footnotesize}
\begin{tabular}{c|ccc}
Box \# & Centre & Eigenvector & Scale\\ \hline
0 & (0.577619,0) & Fixed & 0.004\\
1 & ( -0.527156, 0) & Period Three & 0.0005\\
2 & (-0.524339,-0.00302553)& Period Three & 0.001\\
3 & (-0.458451,-0.0695428) & Period Three &0.002\\
4 & (0.683671,-0.0907786) & Fixed & 0.002\\
5 & (0.565061,0.00974195) & Fixed & 0.002\\
6 & (0.578698,-0.0014086) & Fixed & 0.003
\end{tabular}
\caption{The boxes used to prove the existence of a heteroclinic orbit 
between the period three point and the fixed point, passing through $q_1$. 
Note that the numbers in this table are in the original coordinates. 
The point $q_1$ corresponds to box $\# 3$.
The spanning vectors are the stable and unstable vectors for either the fixed point or the 
period three point. Note, the stable vector is constructed by reflecting 
the unstable vector in $\{u=0\}$. The approximate unstable eigenvectors 
for the fixed and period three points are $(0.788578889012330, 
0.614933602760558)$ and $( 0.682259082166558, -0.731107710380897)$, 
respectively. To construct the corresponding parallelogram the spanning 
vectors are scaled as indicated in the table. The covering sequence is 
$1\Rightarrow1\Rightarrow2\Rightarrow3\Rightarrow4\Rightarrow5 
\Rightarrow6 \Rightarrow 0\Rightarrow 0$.}\label{tHeteroTraj}
\end{footnotesize}
\end{center}

\begin{center}
\begin{footnotesize}
\begin{tabular}{c|ccc}
Box \# & Centre & Eigenvector & Scale\\ \hline
0 & (0.577619,0) & Fixed & 0.002\\
1 & ( -0.527156, 0) & Period Three & 0.00015\\
2 & (-0.527035,-0.000143356) & Period Three & 0.0001(0,9)\\
3 & (-0.523879,-0.00351678) & Period Three & 0.0003\\
4 & (-0.447416,-0.0799474) & Period Three & 0.0003\\
5 & (0.787595,-0.00888005) & Fixed & 0.0005\\
6 & (-0.36687,0.0109682) & Period Three & 0.0003 \\
7 & (0.67407,-0.0819091) & Fixed & 0.001\\
8 & (0.566215,0.00882267) & Fixed & 0.001
\end{tabular}
\caption{The boxes used to prove the existence of a heteroclinic orbit 
between the period three point and the fixed point, passing through $q_2$. 
Note that the numbers in this table are in the original coordinates. 
$k_1^t$ corresponds to box $\# 7$ and $q_2$ to box $\# 4$. The spanning vectors 
are the stable and unstable vectors for either the fixed point or the 
period three point. Note, the stable vector is constructed by reflecting 
the unstable vector in $\{u=0\}$. The approximate unstable eigenvectors 
for the fixed and period three points are $(0.788578889012330, 
0.614933602760558)$ and $( 0.682259082166558, -0.731107710380897)$, 
respectively. To construct the corresponding parallelogram the spanning 
vectors are scaled as indicated in the table. The covering sequence is 
$1\Rightarrow 1\Rightarrow 2\Rightarrow 3\Rightarrow 4\Rightarrow 5 
\Rightarrow 6\Rightarrow 7\Rightarrow 8\Rightarrow 0\Rightarrow 0$.}\label{HeteroTraj2}
\end{footnotesize}
\end{center}
\end{table}
{\flushright {$\Box$}}

Consider a $G\in\bG$, by Lemma $\ref{lDeltaTilde}$ the two components of the Markov partition for $\cK^u$ are contained in the intersections
$$\tilde G(\tD_0) \cap \left( \tilde \cW^s_\loc(\tp_0) \cap \tD_0\right),\quad 
\tilde{G}(\tD_1) \cap \left(\tilde \cW^s_\loc(\tp_0) \cap \tD_0\right).$$

Therefore, our next task will be two find an enclosure of the local stable manifold at the fixed point $p_0$.

By the property of $h$-sets with cone conditions mentioned in Remark \ref{rSMfd}, 
it is possible to propagate stable manifolds backwards through a 
sequence of covering relations. We use this property in the program {\tt encloseStableManifold} to enclose the local 
stable manifold  $\tilde \cW^s_\loc(\tp_0) \cap \tD_0$; estimates of the length of its intersection with $\tilde G_*(\tD)$ are needed to estimate the left and right thicknesses of the Markov  partition for $\cK^u$, and ultimately the Hausdorff dimension.

To enclose the stable manifold of the fixed point, we work in the coordinates $(\ref{linChange})$.

First, we choose $200$ points $q_i$ distributed uniformly on the interval 
$$\left(p_0,p_0+0.00011 (0.788578889012330, -0.614933602760558)\right) \subset \fR^2,$$
and iterate each of them $10$ times with $F_*^{-1}$. These $2200$ points $z_{i,n} \equiv F_*^{-n}(q_i)$, $ 1 \le i\le 200$, $0 \le n \le 10$, provide an approximation of the stable manifold.

We next construct rectangles 
\begin{eqnarray}
\nonumber B_{i,n} \equiv \{(x,u) \in \fR^2: &|&(x-\cP_x \cT(z_{i,n}),u-\cP_u \cT(z_{i,n}))\cdot {\bf e}^u |<0.001,\\
\nonumber  &|& (x-\cP_x \cT(z_{i,n}),u-\cP_u \cT(z_{i,n})) \cdot {\bf e}^s |<0.001 \},
\end{eqnarray}
where
\begin{eqnarray}
\nonumber {\bf e}^u&=&(0.570868623900820281, -0.821041420541974730), \\
\nonumber {\bf e}^s&=&(0.992704972028756258, 0.120568812341277774)
\end{eqnarray}
are the approximate unstable and stable eigenvectors of $D\tilde{F}_*$ at $p_0$, and 
verify that
$$B_{i,10} \cover{\tF} B_{i,9} \cover{\tF} \ldots \cover{\tF} B_{i,0} \cover{\tF}B_0,$$
for all $F\in\bF$, where 
\begin{eqnarray}
\nonumber B_0 \equiv \{(x,u) \in \fR^2: &|&(x-\cP_x \cT(p_0),u-\cP_u \cT(p_0))\cdot {\bf e}^u |<0.001, \\
\nonumber   &|&(x-\cP_x \cT(p_0),u-\cP_u \cT(p_0)) \cdot {\bf e}^s |<0.001 \}.
\end{eqnarray}

To prove that the cone conditions are satisfied we verify the hypothesis of Lemma \ref{lQF} using the routine {\tt checkConeConditions}. This proves that $\left( \bigcup_{i,n} B_{i,n} \right) \cup B_0$ covers the local stable manifold.

The sets $\bk_0^b, \bk_0^t, \bk_1^b,$ and $\bk_1^t$, from Lemma \ref{lKuEnc}, serve as bounds on the end points of the sets $\cK_{0,1}^u$ in the Markov partition $\cP=\{\cK_0^u,\cK_1^u\}$. Therefore,
\begin{eqnarray}
\label{tauL} \tau_L(\cP) \ge {{\rm dist}(\bk_1^b,\bk_1^t) \over \sqrt{1+{\rm Lip}^2} {\rm sup}_{q \in \bk_1^t, q' \in \bk_0^b } |q-q'|} \ge 0.0650166,\\
\label{tauR} \tau_R(\cP) \ge {{\rm dist}(\bk_0^b,\bk_0^t) \over \sqrt{1+{\rm Lip}^2} {\rm sup}_{q \in \bk_1^t, q' \in \bk_0^b } |q-q'|} \ge 0.0514139
\end{eqnarray}
(these quantities are computed  in the routine {\tt Gaps} in \cite{GP}).

Finally, we compute (in routine {\tt Distortion} of \cite{GP}) the quantities $C$, $\gamma$ and $\epsilon$ appearing in the Definition $\ref{classF}$, and verify parts $3)$, $4)$ and $5)$ of this definition. In particular, we have set
$$\epsilon=\sup_{\Delta_0 \cup \Delta_1}\left\{{ \max\{|b|,|c|\} \over |a|-1 } \right\}, \quad C=\sup_{\Delta_0 \cup \Delta_1}\left\{ |a| {\max\{|b|,|c|\} \over |a|-1 } \right\},$$
while $\gamma$ has been chosen as the supremum of all possible values that verify the conditions $(\ref{gamma1})$--$(\ref{gamma2})$ over $\Delta_0 \cup \Delta_1$. To verify part $4)$ we have bound the total variation of $\eta(x,u)=\log|a(x,u)|$ on $\Delta_i'$ from above by ${\rm Area}(\Delta_i') \cdot \sup_{(x,u) \in \Delta_i'} |\nabla \eta(x,u)|.$ 

\begin{lemma}\label{Cge}
Any function $G\in\bG: \Delta_0 \cup \Delta_1 \mapsto \fR^2$ is of class $\cF(C,\gamma,\epsilon)$ with
$$C=16.6, \quad \gamma=47.8, \quad \epsilon=0.88,$$ 
\end{lemma} 

A straightforward implementation of the Duarte Distortion Theorem and Lemma $\ref{thicknesslemma}$ gives us a bound on the Hausdorff dimension of the unstable Cantor set:

\begin{corollary}\label{cHdfDim}
For all $G \in \bG$, the distortion of the unstable Cantor set $\cK^u_G$ is less or equal to
$$D(C,\gamma,\epsilon)=3749.3,$$
while its Hausdorff dimension satisfies:
$${\rm dim}_H(\cK^u_G) \ge \left\{ {\log\left(1+{e^{-D} \tau_R(\cP) \over 1+e^{D} \tau_L(\cP)} \right) \over \log\left(1+{1+e^{D} \tau_R(\cP) \over e^{-D} \tau_L(\cP)} \right) },   {\log\left(1+{e^{-D} \tau_L(\cP) \over 1+e^{D} \tau_R(\cP)} \right) \over \log\left(1+{1+e^{D} \tau_L(\cP) \over e^{-D} \tau_R(\cP)} \right) }    \right\} \approx {1 \over 2 D} e^{-2 D}{\tau_R(\cP) \over \tau_L(\cP)},$$
where $\tau_L(\cP)$ and $\tau_R$ are as in $(\ref{tauL})$--$(\ref{tauR})$.

The Hausdorff dimension of the set $\cC_G$ satisfies 
$${\rm dim}_H(\cC_G) \ge 2 {\rm dim}_H(\cK^u_G) \approx 0.00013 \, e^{-7499}.$$

\end{corollary}
\noindent {\it Proof.}
The fact that the stable and unstable foliations are transversal (see Lemma $\ref{transversal_cones}$) implies that $\cC_G = \cK^u_G \times \cK^s_G$. Reversibility and the symmetry of the Markov partition $\Delta_0 \cap \Delta_1$ implies that $\cK^s_G=T(\cK^u_G)$, and ${\rm dim}_H(\cK^s_G)={\rm dim}_H(\cK^u_G)$, which in turn implies the last claim.
\newline
{\flushright $\Box$}

Corollary \ref{cHdfDim} together with Theorem \ref{tP3TopHorse} proves the first part of the second Main Theorem.

Let us now construct a convergent sequence of approximations of the Cantor sets $\cC_G$.
As before, we denote 
$\Delta_0$ and $\Delta_1$ the two full components of the Markov partition as in Definition $\ref{classF}$ and Lemma $\ref{lDeltaTilde}$, $\Delta \equiv \Delta_0 \cup \Delta_1$. Define recursively:
$$ \cU_G^1 \equiv G(\Delta) \cap G^{-1}(\Delta) \quad {\rm and} \quad \cU_G^k \equiv  
G(\cU_G^{k-1}) \cap G^{-1} (\cU_G^{k-1}).$$  

Each of the sets $\cU_G^k$ contains $2 \cdot 4^k$ components 
$\cU_G^{k,n}$, $n=1..2 \cdot 4^k$.

Recall that $\partial \Delta_j=b_j \cup t_j \cup l_j \cup r_j$, $j=0,1$. By definition of $\cU_G^{k,n}$, the boundary $\partial \cU_G^{k,n}$ consists of two images, $t^{k,n}$ and $b^{k,n}$, of some subsets of $t_j$ and $b_j$, $j=0,1$, under $G^k$, and two images, $l^{k,n}$ and $r^{k,n}$, of subsets of $l_j$ and $r_j$, $j=0,1$, under $G^{-k}$. We will refer to $t^{k,n}$, $b^{k,n}$,  $l^{k,n}$ and  $r^{k,n}$ as  ``edges'' or ``sides'', and to their intersections as ``corners'' of $\cU_G^{k,n}$.

\begin{lemma}\label{lMSBds}
Let 
$$\rho_{k,n}={\rm sup}_{B_\rho \subset \cU_G^{k,n}} (\rho),$$
and set $\rho_k=\min_n \{ \rho_{k,n} \}$. There exist constants $C>0$ and $c>0$ such that
\begin{eqnarray}
\nonumber {\rm diam}\left(\cU_G^{k,n}\right)&\le& C \, A_-^{-k} \equiv C \kappa_+^k,\\
\nonumber \phantom{aaaaaaaa} \rho_k &\ge& c \, A_+^{-k} \equiv c \kappa_-^k ,
\end{eqnarray} 
where $A_\pm$ are as in $(\ref{Apm})$.
\end{lemma}
\noindent {\it Proof.} The length of an edge of $\cU_G^{k,n}$ is bounded by $c A^{-k}$ for some $A \in \bA$. This follows from the fact that $\Delta_j$, $j=0,1$, are constructed so that the tangent vectors to $t_j$ and $b_j$, $j=0,1$, are contained in $\cC_u$, and those to $l_j$ and $r_j$, $j=0,1$, --- in $\cC_s$ (see $(\ref{Cu})$ and $(\ref{Cs})$), and from the fact that an edge of $\cU_G^{k,n}$ is contained in the image of an edge of  $\cU_G^{k-1,j}$ for some $j$. 

Consider an edge of $\cU_G^k$, say $t^{k,n}$. Set
$$
 \alpha_t^+ \equiv \sup_{p \in t^{k,n}} \left\{ {\rm maximum \,\,} \angle{\bf v}: {\bf v} \in \cC_u(p)  \right\}, \quad  \alpha_t^- \equiv \inf_{p \in t^{k,n}} \left\{ {\rm minimum\,\,} \angle{\bf v}: {\bf v} \in \cC_u(p)  \right\},
$$
where $\angle {\bf v}$ signifies the angle between vector ${\bf v}$ and the x-axis, measured, say, counterclockwise. Similarly for $b^{k,n}$, $l^{k,n}$ and $r^{k,n}$. Straightforward geometric considerations 
demonstrate that the angles $\alpha^\pm_{t,b,l,r}$ specify two quadrilaterals, containing and contained in $\cU_G^{k,n}$, respectively. 

Since the angle between any ${\bf v} \in \cC^u(p)$ and  any ${\bf u} \in \cC^s(p)$ is bounded from below (see Lemma $\ref{transversal_cones}$), the angles of the quadrilaterals are bounded away from $0$ and $\pi$, and their diameters satisfy:
\begin{eqnarray}
\nonumber c \, A_+^{-k} &\leq& {\rm 
diam}\left({\rm Inner \,\, Quadrilateral}\right),\\
\nonumber C \, A_-^{-k} &\geq& {\rm 
diam}\left({\rm Outer \,\, Quadrilateral}\right).
\end{eqnarray}

The claim follows.
\newline
{\flushright $\Box$}

\begin{corollary}\label{cHdfDimU}
The Hausdorff dimension of the Cantor set $\cC_G$ satisfies 
$${\rm dim}_H(\cC_G) \leq \frac{\log(4)}{-\log{\kappa_+}} \le 0.7673.$$
\end{corollary}
\noindent {\it Proof.}
Clearly $\cC_G\subset \cU^k_G$, for all $k$. As above, $\cU_G^k=\bigcup_{n=1}^{2 \cdot 4^k} \cU_G^{k,n}$, where 
$$
{\rm diam}(\cU_G^{k,n}) \leq {\rm const} \, \kappa_+^k, \quad \textrm{ for all  }1\leq n\leq 2 \cdot 4^k.
$$
Hence, for any $k>1$,
$C_d^H(\cC_G)\leq {\rm const} \, 4^k \kappa_+^{kd}$, and if $4 \kappa_+^d<1$, then $C_d^H(\cC_G)=0$.   
\newline
{\flushright $\Box$}

\section{Overview of the programs}
We will now give a brief overview of the programs used in the proofs.

A subset of the domain of a map $F\in \bF$ is determined using the program {\tt findDomain} of \cite{JP}, used in Lemma \ref{lfDomRe}. {\tt findDomain} discretizes the domain of the generating function in the $x$-domain, and computes rigorous bounds on the maximum and minimum of the generating function, $s(y,x)$, on this slice of the domain. During the process of computing the maximum and minimum, the program also verifies that the conditions of the implicit function theorem are verified on the slice.

To verify that a pair of $h$-sets satisfies the hypothesis of Theorem \ref{tCovRel}, we use the routine {\tt checkCoveringRelations} of \cite{JP}. It discretizes the boundary of a $h$-set, and verifies that each piece of the boundary is mapped so that the covering relations hold. 

To verify that a pair of $h$-sets with cones satisfies the hypothesis of Lemma \ref{lQF}, we use the routine {\tt checkConeConditions} of \cite{JP}. It computes the derivative of the map on the initial $h$-set in appropriate coordinates, and checks the positive definite condition of the quadratic form $V$ in the statement of the lemma using the Sylvester criterion.

The sequences of $h$-sets used to enclose the stable manifold in Section \ref{HDC}, are constructed in the program {\tt encloseStableManifold} of \cite{JP}. This program uses the routines {\tt checkCoveringRelations} and {\tt checkConeConditions}, to verify that on each sequence, the hypotheses of Theorem \ref{tCovRel} and Lemma \ref{lQF}, respectively, are satisfied.

As we have already mentioned parts of Theorem $\ref{tP3TopHorse}$, Lemmas $\ref{transversal_cones}$ and  $\ref{Cge}$ are proved with the help of the routine {\tt Distortion} of \cite{GP}. The input of the routine is a set of parameters that define the components $\Delta_0'$ and $\Delta_1'$ (see Table \ref{bigHorseshoe}) and the cone fields (see  $(\ref{tCu1})$--$(\ref{tCu2})$). The routine verifies that the cone fields are invariant and computes the expansion rates $(\ref{Apm})$. It also computes the quantities $C$, $\gamma$ and $\epsilon$ appearing in Lemma $\ref{Cge}$.

The thicknesses $(\ref{tauL})$ and $(\ref{tauR})$ of the Markov partition for the horseshoe of $\bG$ are computed in the routine {\tt Gaps} of \cite{GP}.

\section{Concluding remarks}

An obvious question is whether the Hausdorff dimension of the set $\cC_G$ for all $G$ in some subset of $\bG$ is independent of $G$. 

In the satellite work \cite{GJ} we were able to use a renormalization approach to show that the Hausdorff dimension of the ``stable'' set, that is the set on which the maximal Lyapunov exponent is zero, is indeed invariant for a subset of {\it infinitely renormalizable} maps. An essential ingredient of that proof is the fact that the renormalizations of all infinitely renormalizable maps in a neighbourhood of the fixed point converge to that fixed point on their domain of analyticity $\cD$. 

In a somewhat similar fashion, the invariance of the Hausdorff dimension of the set $\cC_G$ could be demonstrated if one can show that renormalizations of the infinitely renormalizable maps converge to $G_*$ on neighbourhoods $\cV_k$ of the rescalings $\Lambda^{-1}_{k,G}(\cC_G)$ of the hyperbolic sets. However, $\cD \cup_{k=1..\infty} \cV_k$ is not a connected domain, and convergence on $\cD$ does not imply that  on $\cup_{k=1..\infty} \cV_k$.

\ack
The authors would like to thank Daniel Wilczak for many interesting 
discussions and his help with the usage of the software package 
\cite{CAPD}. We are also grateful to Hans Koch for his many useful 
insights into period-doubling in area-preserving maps, as well as for 
his great help with understanding the original computer assisted proof 
\cite{EKW2}.


\end{document}